\documentclass[article,3p]{elsarticle} 

\usepackage{amsmath,amsfonts,amssymb,amsbsy,amscd}

\usepackage{mathrsfs,array}
\usepackage{verbatim}
\usepackage{subcaption}
\usepackage{amsthm}
\usepackage{tikz}

\makeatletter
\def\ps@pprintTitle{%
 \let\@oddhead\@empty
 \let\@evenhead\@empty
 \def\@oddfoot{}%
 \let\@evenfoot\@oddfoot}
\makeatother

\newdefinition{rmk}{Remark}
\newproof{pf}{Proof}
\newproof{pot}{Proof of Theorem \ref{thm2}}

\def\B#1{\mbox{\boldmath{$#1$}}}

\newcommand{\pd}{\partial}

\newcommand{\bx}{\B{x}}

\newcommand{\onehalf}{\mbox{$\frac{1}{2}$}}
\newcommand{\VV}{\mathcal{V}}
\newcommand{\WW}{\mathcal{W}}

\def\be{\begin{equation}}
\def\ee{\end{equation}}
\def\ba{\begin{array}}
\def\ea{\end{array}}
\def\bea{\begin{eqnarray}}
\def\eea{\end{eqnarray}}
\def\beas{\begin{eqnarray*}}
\def\eeas{\end{eqnarray*}}
\newcommand{\bseq}{\begin{subequations}}
\newcommand{\eseq}{\end{subequations}}

\biboptions{numbers,sort&compress}

\usepackage{hyperref}
\hypersetup{
    colorlinks,
    citecolor=ck,
    filecolor=ck,
    linkcolor=ck,
    urlcolor=ck
}

\usepackage[final]{changes}

\definechangesauthor[color=blue]{Rev.1}
\definechangesauthor[color=orange]{Rev.2}
\definechangesauthor[color=red]{Both Rev}
\definechangesauthor[color=purple]{Authors}

\journal{Computer Methods in Applied Mechanics and Engineering}

\begin{document}

\begin{frontmatter}

\title{{\large{Correct energy evolution of stabilized formulations: 
The relation between VMS, SUPG and GLS via dynamic orthogonal 
small-scales and isogeometric analysis.}}\\ {\large{II: The incompressible Navier--Stokes equations}}}

\author{M.F.P. ten Eikelder\corref{cor1}}
\ead{m.f.p.teneikelder@tudelft.nl}
\cortext[cor1]{Corresponding author}
\author{I. Akkerman}
\ead{i.akkerman@tudelft.nl}

\address{Delft University of Technology, Department of Mechanical, 
Maritime and Materials Engineering, P.O. Box 5, 2600 AA Delft, The Netherlands}
\date{}

\begin{abstract}
This paper presents the construction of a correct-energy stabilized finite 
element method for the incompressible Navier-Stokes equations. 
The framework of the methodology and the correct-energy concept have
been developed in the convective--diffusive context in the preceding 
paper [M.F.P. ten Eikelder, I. Akkerman, Correct energy evolution of
stabilized formulations: The relation between VMS, SUPG and GLS
via dynamic orthogonal small-scales and isogeometric analysis. 
I: The convective--diffusive context, Comput. Methods Appl. Mech. Engrg. 331 (2018) 259--280]. 
The current work extends ideas of the preceding paper to build a stabilized method within the
variational multiscale (VMS) setting which displays correct-energy behavior. 
Similar to the convection--diffusion case, a key ingredient is the proper 
dynamic and orthogonal behavior of the small-scales. This is demanded 
for correct energy behavior and links the VMS framework to the 
streamline-upwind Petrov-Galerkin (SUPG) and the 
Galerkin/least-squares method (GLS). 

The presented method is a Galerkin/least-squares formulation with 
dynamic divergence-free small-scales (GLSDD). It is locally 
mass-conservative for both the large- and small-scales separately.
 In addition, it locally conserves linear and angular momentum. 
 The computations require and employ NURBS-based isogeometric 
 analysis for the spatial discretization. The resulting formulation 
 numerically shows improved energy behavior for turbulent flows 
 comparing with the original VMS method.
\end{abstract}


\end{frontmatter}

\section{Introduction}
\label{sec:Introduction}

The creation of artificial energy in numerical methods is undesirable from 
both a physical and a numerical stability point of view. Therefore methods 
precluding this deficiency are often sought after. This work continues the 
construction of the correct-energy displaying stabilized finite element methods. 
The first episode \cite{EAk17} exposes the developed methodology in the 
convective--diffusive context. The current study deals with the incompressible 
Navier--Stokes equations and is the second piece of work within the framework.
 The setup of this paper is closely related to that of \cite{EAk17}. In particular, 
 the \textit{correct-energy} demand is the same, thus it represents that the 
 method (i) does not create artificial energy and (ii) closely resembles the 
 energy evolution of the continuous setting. The precise definition is stated 
 in Section \ref{sec:Towards correct energy behavior}. 
What sets the Navier--Stokes problem apart from convection--diffusion case 
is the inclusion of the incompressibility constraint. 
In this work we use a divergence-conforming basis which allows exact 
pointwise satisfaction of this constraint.
This is considered a beneficial property. Therefore it is added as a design criterion.
\added[id=Rev.1]{In a two-phase context this property is essential for correct energy behavior \cite{akkerman2018toward}.}

\subsection{Contributions of this work}
This paper derives a novel VMS formulation which exhibits the correct 
energy behavior and to this purpose combines several ingredients. 
The final formulation is summarized in \ref{Appendix: Galerkin/least-squares 
formulation with dynamic divergence-free small-scales}. The new method 
is a residual-based approach that employs (i) dynamic behavior of the 
small-scales, (ii) solenoidal NURBS basis functions and (iii) a Lagrange-multiplier 
construction to ensure the incompressibility of the small-scale velocities. 
The formulation is of skew-symmetric type, rather than conservative, which 
is motivated by both the correct-energy demand and its improved behavior 
in the single scale setting (i.e. the Galerkin method) \cite{HuWe05}. 
Moreover, \added[id=Rev.1]{the formulation reduces to a Galerkin formulation in case of a vanishing Reynolds number due to a Stokes-projector.}
The use of dynamic small-scales, firstly proposed in \cite{Cod02}, is also driven
 from an energy point of view. In addition, it leads to global momentum 
 conservation and \added[id=Rev.1]{the numerical results of \cite{CodPriGuaBad07} show improved behavior of the dynamic small-scales} with respect 
 to their static counterpart.

\subsection{Context}
\added[id=Rev.1]{This work falls within the variational multiscale framework \cite{Hug95, Hug98}. 
The basic idea of this method is to split solution into the large/resolved-scales and small/unresolved-scales.
The small-scales are modeled in terms of (the residual of) the large scales and substituted into the equation for the large-scales.
This approach was first applied in a residual-based LES context to incompressible 
turbulence computations in \cite{BaCaCoHu07}. 
The VMS methodology has enjoyed a lot of progress since then. For an overview of the development 
consult the review paper \cite{codina2017variational}.

Our work is not the first to analyze the energy behavior of the VMS method.
A spectral analysis of the VMS method can be found in \cite{wang2010spectral}.
That paper proves dissipation of the model terms under restrictive conditions. 
Additional to the optimality projector, they require $L_2$-orthogonality of the large- and small-scales. 
This condition naturally leads to the use of spectral methods.

Principe et al. \cite{principe2010dissipative} provide a precise definition 
of the numerical dissipation within the variational multiscale context for 
incompressible flows. Equally important, they numerically show that the 
concept of dynamic small-scales, which we apply in this work, is able 
to model turbulence.

Colom\'{e}s et al. \cite{colomes2015assessment} assess the performance of several VMS methods for turbulent flow problems and provide an energy analysis of these methods.
They conclude that algebraic subgrid scales (ASGS) and orthogonal subscales (OSS) yield similar results, whereas the latter one is more convenient in terms of numerical performance.

We build onto \cite{wang2010spectral,principe2010dissipative,colomes2015assessment} without requiring $L_2$-orthogonality.
 Therefore we are not restricted to the use of spectral methods, while retaining a strict energy relation.

Other} recent related work includes the IGA divergence-conforming VMS method 
of Opstal et al. in \cite{OpYanCoEvKvBa17}. They also employ an 
$H_0^1$-orthogonality between the velocity large- and small-scales \added[id=Rev.1]{on a local level}.
 Our work deviates from 
\cite{OpYanCoEvKvBa17} in that we motivate the required orthogonalities 
with the correct energy demand. Furthermore, our work distinguishes 
itself by enforcing the divergence-free velocity small-scales with a 
Lagrange-multiplier construction. We believe the Stokes orthogonality 
between the large- and small-scales is a natural path to take, since it 
reduces the scheme to the Galerkin method in the vanishing Reynolds number limit. 

The discretizations throughout this work are based on the isogeometric 
analysis (IGA) concept, proposed by Hughes et al. in \cite{HuCoBa04}. 
This idea integrates the historically distinct fields of computer aided 
design (CAD) and finite element analysis. Isogeometric analysis rapidly 
became a valuable tool in computational fluid dynamics, in particular in 
turbulence computations. It provides several advantages over standard 
finite element analysis, including an exact description of CAD geometries, 
increased robustness and superior approximation properties 
\cite{HuCoBa04, LEBEH09, BACHH07}. \added[id=Both Rev]{This work requires in particular 
inf--sup stable discretizations for which we use \cite{Evans13steadyNS, Evans13unsteadyNS}.
Moreover these spaces allow the pointwise satisfaction of the incompressibility constraint.}
The smooth NURBS basis functions are convenient for the computation of second derivatives. 

\subsection{Outline}
The organization of this paper in Section \ref{sec:GE} and \ref{sec:EESSM} 
is very comparable with that of the convective--diffusive context \cite{EAk17}, 
and at some points mirrors it. The purpose thereof is (i) to indicate the great 
similarities of the methodologies and (ii) to clarify the approach. 
The remainder of this paper presents the actual construction of a stabilized 
variational formulation for the incompressible Navier--Stokes equations 
which displays correct-energy behavior. We summarize it as follows. 
Section \ref{sec:GE} states the continuous form of the governing 
incompressible flow equations, both in the strong formulation and the 
standard weak formulation. It additionally provides the energy evolution 
of the continuous equation, in both global and local form. 
Section \ref{sec:EESSM} discusses the energy evolution of the variational 
multiscale approach with dynamic small-scales. The path toward correct 
energy behavior actually starts in Section \ref{sec:Towards correct energy behavior}. 
This Section presents the required orthogonality of the large-scales and small-scales. 
This converts the residual-based variational multiscale method into the 
Galerkin/least-squares method with the correct energy behavior. 
Section \ref{sec:CP} presents conservation properties of the method. 
Section \ref{sec:ns_case} provides a computational test case, namely a 
three-dimensional Taylor--Green vortical flow.
In particular it examines the energy behavior and compares the novel 
method with the standard VMS method with static small-scales \cite{BaCaCoHu07}. 
The calculations employ the generalized-$\alpha$ method with favorable 
energy behavior which is also discussed in \cite{EAk17}. 
In Section \ref{sec:CONC}, we wrap up and present avenues for future research.

\section{The continuous incompressible Navier--Stokes equations}
\label{sec:GE}

\subsection{Strong formulation}
\label{sec:GE, subsec:SF}

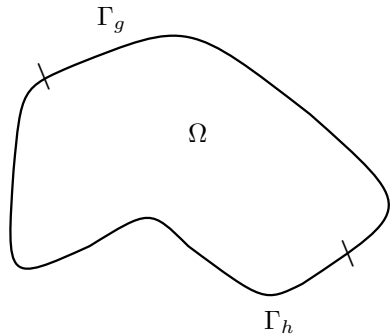
\begin{figure}[h]
  \begin{center}
\begin{tikzpicture}
\draw[line width=0.3mm, black ] (2,0.4) .. controls (3.5,1) .. (5.1, -0.25) .. controls (6.5,-1.5) .. (5, -2.5) .. controls (4.5,-2.75) .. (3.5, -2).. controls (3.0,-1.5) .. (2.2, -2) .. controls (1.1,-2.5) .. (1.2, -1.0) .. controls (1.3,0.1) .. (2.0,0.4);
\node[text width=3cm] at (5.0,-0.5) {$\Omega$};
\node[text width=3cm] at (6.0,-3.0) {$\Gamma_h$};
\node[text width=3cm] at (3.8,1.0) {$\Gamma_g$};
\node[text width=3cm] at (7.0,-2.1) {$\boldsymbol{\backslash}$};
\node[text width=3cm] at (3.0,0.25) {$\boldsymbol{\backslash}$};
\end{tikzpicture}
  \end{center}
  \caption{Spatial domain $\Omega$ with its boundaries 
  $\Gamma = \Gamma_g \cup \Gamma_h$. This is the same figure as in \cite{EAk17}.}
  \label{fig:domain1}
\end{figure}

Let $\Omega \in \mathbb{R}^d$, \added[id=Rev.2]{$d=2,3$}, denote the spatial domain 
and $\partial \Omega = \Gamma=\Gamma_g \cup \Gamma_h$ its 
boundary, see Figure \ref{fig:domain1}.
\\
\\
The problem consists of solving the incompressible Navier--Stokes equations 
governing the fluid flow, which read in strong form
\begin{subequations}
  \label{sec:GE, subsec:SF, NS Strong}
  \begin{alignat}{1}
     \pd_t  \mathbf{u} +  \nabla \cdot \left( \mathbf{u} \otimes \mathbf{u}\right) 
     + \nabla p- \nabla \cdot \left( 2 \nu \nabla^s \mathbf{u} \right)= \mathbf{f} & \quad \text{in} \quad \Omega\times \mathcal{I}, 
     \label{sec:GE, subsec:SF, NS Strong, mom eq}\\
    \nabla \cdot \mathbf{u} = 0 & \quad \text{in} \quad \Omega \times \mathcal{I}, \label{sec:GE, subsec:SF, NS Strong, div free eq}\\
    \mathbf{u} =\mathbf{g} & \quad \text{in} \quad  \Gamma_g \times \mathcal{I}, \label{sec:GE, subsec:SF, NS Strong, dir bc}\\
    -u_n^-\mathbf{u}-p\mathbf{n}+\nu\pd_n \mathbf{u} =\mathbf{h} & \quad \text{in} \quad \Gamma_h \times \mathcal{I}, 
    \label{sec:GE, subsec:SF, NS Strong, neum bc}\\
    \mathbf{u}(\mathbf{x},0) = \mathbf{u}_0(\mathbf{x})  & \quad \text{in} \quad \Omega, \label{sec:GE, subsec:SF, NS Strong, IC}
  \end{alignat}
\end{subequations}
for the velocity $\mathbf{u}: \Omega \times \mathcal{I} \rightarrow \mathbb{R}^d$ 
and the pressure divided by the density $p: \Omega \times \mathcal{I} \rightarrow \mathbb{R}$. 
A constant density is assumed. Eqs. (\ref{sec:GE, subsec:SF, NS Strong, mom eq})-(\ref{sec:GE, subsec:SF, NS Strong, IC}) 
describe the balance of linear momentum, the conservation of mass, the inhomogeneous Dirichlet 
boundary condition, the traction boundary condition and the initial conditions, respectively. The 
spatial coordinate denotes $\mathbf{x} \in \Omega$ and the time denotes $t \in \mathcal{I}=(0,T)$ 
with end time $T>0$. The given dynamic viscosity is $\nu: \Omega \rightarrow \mathbb{R}^+$, the 
body force is $\mathbf{f}: \Omega \times \mathcal{I} \rightarrow \mathbb{R}^d$, the initial velocity 
is $\mathbf{u}_0: \Omega \rightarrow \mathbb{R}^d$ and the boundary data are 
$\mathbf{g}: \Gamma_g \times \mathcal{I} \rightarrow \mathbb{R}^d$ and 
$\mathbf{h}: \Gamma_h \times \mathcal{I} \rightarrow \mathbb{R}^d$. 
We assume a zero-average pressure for all $t \in \mathcal{I}$ \added[id=Both Rev]{in case of an empty Neumann boundary}. The normal velocity 
denotes $u_n=\mathbf{u}\cdot\mathbf{n}$ with positive and negative parts 
$u_n^{\pm}=\tfrac{1}{2}(u_n\pm|u_n|)$. The various derivative operators are the 
temporal one $\pd_t$, the symmetric gradient $\nabla^s\cdot=\tfrac{1}{2}\left(\nabla\cdot + \nabla^T\cdot\right)$ 
and the normal gradient $\pd_n =\mathbf{n}\cdot \nabla $, with $\mathbf{n}$ the outward unit normal.
\subsection{Weak formulation}
\label{sec:GE, subsec:SWF}

Let $\WW^0$ denote the trial weighting function space satisfying the homogeneous 
Dirichlet conditions on $\mathbf{u}$ and $\WW^g$ the trial solution space with 
non-homogeneous Dirichlet conditions on $\mathbf{u}$. The standard variational formulation writes:\\

\textit{Find $\left\{\mathbf{u}, p\right\} \in \WW^g$ such that for all $\left\{\mathbf{w}, q\right\} \in \WW^0$,}
\begin{subequations}
  \label{sec:GE, subsec:SWF, standard weak form1}
  \begin{alignat}{2}
     B_{\Omega,\Gamma_h}\left(\left\{\mathbf{u}, p\right\},\left\{\mathbf{w}, q\right\}\right)=& \added[id=Rev.1]{L_{\Omega,\Gamma_h}}\left(\left\{\mathbf{w}, q\right\}\right),
       \end{alignat}
   \indent \textit{where}
       \begin{alignat}{2}
     B_{D,\Gamma_h}\left(\left\{\mathbf{u}, p\right\},\left\{\mathbf{w}, q\right\}\right)=&B_{D}\left(\left\{\mathbf{u}, p\right\},\left\{\mathbf{w}, q\right\}\right)+\left(\mathbf{w},u_n^+\mathbf{u}\right)_{\Gamma_h\left(D\right)},\\
     \added[id=Rev.1]{L_{D,\Gamma_h}}\left(\left\{\mathbf{w}, q\right\}\right)=&\added[id=Rev.1]{L_{D}}\left(\left\{\mathbf{w}, q\right\}\right) +\left(\mathbf{w},\mathbf{h}\right)_{\Gamma_h\left(D\right)},\\
     B_D\left(\left\{\mathbf{u}, p\right\},\left\{\mathbf{w}, q\right\}\right)=&\left( \mathbf{w}, \pd_t  \mathbf{u}  \right)_{D}
     -(\nabla \mathbf{w}, \mathbf{u} \otimes \mathbf{u})_{D}+(\nabla \mathbf{w}, 2 \nu \nabla^s \mathbf{u} )_{D}\nonumber\\
     &+(q,\nabla \cdot \mathbf{u})_{D}- (\nabla \cdot \mathbf{w}, p )_{D},\\
     \added[id=Rev.1]{L_D}\left(\left\{\mathbf{w}, q\right\}\right)=&(\mathbf{w},\mathbf{f})_{D}.
  \end{alignat}
\end{subequations}
Here $B_D$ is the bilinear form and $\left(\cdot,\cdot\right)_D$ is the 
$L^2\left(D\right)$ inner product over $D$. The Dirichlet and traction 
boundary of domain $D$ denote
 $\Gamma_g(D):=\Gamma_g\cap\partial D$ 
 and $\Gamma_h(D):=\Gamma_h\cap\partial D$ 
 respectively. The strong (\ref{sec:GE, subsec:SF, NS Strong}) and the 
 weak formulation (\ref{sec:GE, subsec:SWF, standard weak form1}) 
 are equivalent for smooth solutions.

\subsubsection*{Remark}
The variational form (\ref{sec:GE, subsec:SWF, standard weak form1}) 
is of conservative type: the incompressibility constraint 
(\ref{sec:GE, subsec:SF, NS Strong, div free eq}) is not directly 
employed in the convective terms. A discretization of the conservative
 form may lead to spurious oscillations caused by the error in the 
 incompressibility constraint acting as a distribution of sinks and sources.
  Employing (\ref{sec:GE, subsec:SF, NS Strong, div free eq}) can be 
  used to generate a \textit{convective form} which is sometimes 
  preferred and often adopted in Galerkin computations \cite{HuWe05}. 
  Here we write the variational formulation of \textit{skew-symmetric} 
  type which will be used in Section \ref{sec:Towards correct energy behavior}:\\

\textit{Find $\left\{\mathbf{u}, p\right\} \in \WW^g$ 
such that for all $\left\{\mathbf{w}, q\right\} \in \WW^0$,}
\begin{subequations}
  \label{sec:GE, subsec:SWF, standard weak form1: convective}
  \begin{alignat}{2}
     C_{\Omega,\Gamma_h}\left(\left\{\mathbf{u}, p\right\},\left\{\mathbf{w}, q\right\}\right)=& \added[id=Rev.1]{L_{\Omega,\Gamma_h}}\left(\left\{\mathbf{w}, q\right\}\right),
\end{alignat}
\indent \textit{where}
\begin{alignat}{2}
     C_{D,\Gamma_h}\left(\left\{\mathbf{u}, p\right\},\left\{\mathbf{w}, q\right\}\right)=&C_{D}\left(\left\{\mathbf{u}, p\right\},\left\{\mathbf{w}, q\right\}\right)+\tfrac{1}{2}\left(\mathbf{w},|u_n|\mathbf{u}\right)_{\Gamma_h\left(D\right)},\\
     C_{D}\left(\left\{\mathbf{u}, p\right\},\left\{\mathbf{w}, q\right\}\right)=&\left( \mathbf{w}, \pd_t  \mathbf{u}  \right)_{D}
     +\tfrac{1}{2}(\mathbf{w}, \mathbf{u} \cdot \nabla \mathbf{u})_{D}-\tfrac{1}{2}(\mathbf{u}\cdot\nabla\mathbf{w}, \mathbf{u})_{D}
     +(\nabla \mathbf{w}, 2 \nu \nabla^s \mathbf{u} )_{D}\nonumber\\
     &+(q,\nabla \cdot \mathbf{u})_{D}- (\nabla \cdot \mathbf{w}, p )_{D}.
  \end{alignat}
\end{subequations}
Again, this form is equivalent to the strong form (\ref{sec:GE, subsec:SF, NS Strong}). 
Form (\ref{sec:GE, subsec:SWF, standard weak form1: convective}) does not 
possess all conservation properties when discretized in a standard way. 
However, this can be restored using a multiscale split, see \cite{HuWe05} 
for details. In the following we continue with the conservative form 
(\ref{sec:GE, subsec:SWF, standard weak form1}).\\

To obtain the energy evolution linked to (\ref{sec:GE, subsec:SF, NS Strong}) 
we want to substitute $\mathbf{w}=\mathbf{u}$. This is not possible in 
(\ref{sec:GE, subsec:SWF, standard weak form1}) due to the different 
boundary conditions of the solution and test function spaces. 
The enforcement of the Dirichlet boundary conditions in the spaces 
bypasses when employing a Lagrange multiplier construction. 
This converts the variational formulation into a \textit{mixed formulation}:\\

\textit{Find $\left(\left\{\mathbf{u}, p\right\}, \boldsymbol{\lambda}_\Omega\right) \in \WW \times \VV$ 
such that for all $\left(\left\{\mathbf{w}, q\right\}, \boldsymbol{\vartheta} \right) \in \WW \times \VV$,}
\begin{align}\label{sec:GE, subsec:SWF, standard weak form1 mixed}
\left(\boldsymbol{\lambda}_\Omega,\mathbf{w}\right)_{\Gamma_g} = B_{\Omega,\Gamma_h}\left(\left\{\mathbf{u}, p \right\},\left\{\mathbf{w}, q \right\}\right)
-\added[id=Rev.1]{L_{\Omega,\Gamma_h}}(\left\{\mathbf{w},q\right\})+\left(\boldsymbol{\vartheta},\mathbf{u}-\mathbf{g}\right)_{\Gamma_g}.
\end{align}
Here $\WW$ is the unrestricted space used for the solution and test 
functions and $\VV$ is a suitable Lagrange multiplier space. 
Section \ref{sec:GE, subsec:GEE} employs formulation 
(\ref{sec:GE, subsec:SWF, standard weak form1 mixed}) to derive the 
corresponding global energy statement. The equivalence of this form 
with the strong form (\ref{sec:GE, subsec:SF, NS Strong}) follows from 
Green's formula and an appropriate choice of the weighting functions. 
The expression of the Lagrange multiplier is a by-product of this execution and yields
\begin{align}\label{sec:GE, subsec:SWF, LM}
  \boldsymbol{\lambda}_\Omega =-\added[id=Authors]{\tfrac{1}{2}}u_n \mathbf{u} - p \mathbf{n} + \nu \pd_n \mathbf{u}.
\end{align}
The multiplier can be interpreted as an auxiliary flux with a convective, 
a pressure and a viscous contribution. Consult \cite{HEML00} for details 
about auxiliary fluxes in weak formulations.
\subsubsection*{Remark}
\added[id=Authors]{Note that we get the same expression when employing the skew-symmetric form (\ref{sec:GE, subsec:SWF, standard weak form1: convective}).}

\subsection{Global energy evolution}\label{sec:GE, subsec:GEE}
The evolution of the global energy follows when substituting 
$\left(\left\{\mathbf{w},q\right\},\boldsymbol{\vartheta}\right)=\left(\left\{\mathbf{u},p\right\},\boldsymbol{\lambda}_{\Omega}\right)$ 
in (\ref{sec:GE, subsec:SWF, standard weak form1 mixed}). 
Employing Green's formula and the strong incompressibility constraint 
(\ref{sec:GE, subsec:SF, NS Strong, div free eq}) we see that the 
convective term only contributes to the energy evolution via a boundary term. 
The global energy, which is defined as 
$E_{\Omega}:=\tfrac{1}{2}(\mathbf{u}, \mathbf{u})_{\Omega}$, evolves as
\begin{align}\label{sec:GE, subsec:EE, energy1}
  \dfrac{{\rm d}}{{\rm d}t} E_{\Omega}=-\|\nu^{1/2}\nabla \mathbf{u}\|^2_{\Omega}+ \left(\mathbf{u},\mathbf{f}\right)_{\Omega}-(1,F_\Omega)_{\Gamma},
\end{align}
where $ \dfrac{{\rm d}}{{\rm d}t}$ is the time derivative and 
$\|\cdot\|^2_{D}$ defines the standard $L^2$-norm over $D$. 
The flux reads:
\begin{align}
  \label{eq:glb_flux}
F_\Omega =
\left \{ \begin{array}{lc}
  -  \mathbf{g}\cdot \boldsymbol{\lambda}_\Omega &\text{on}~~~ \Gamma_g,  \\
  |u_n|e  - \mathbf{u}\cdot\mathbf{h}&\text{on} ~~~ \Gamma_h,
  \end{array} \right .
\end{align}
with $e:=\frac{1}{2}\mathbf{u}\cdot\mathbf{u}$ the pointwise energy.
 The terms of (\ref{sec:GE, subsec:EE, energy1}) represent from 
 left to right: (i) the energy loss due to viscous molecular dissipation, 
 (ii) the power exerted by the body force and (iii) the energy change 
 due to the boundary conditions. Substitution of the Lagrange 
 multiplier and the boundary conditions leads to the expected 
 expression of the flux
\begin{align}
  \label{eq:glb_flux alt}
F_\Omega =  u_n (e+p) - \nu \pd_n e \quad \text{on} ~~~ \Gamma.
\end{align}
These terms represent the convective and viscous flux as well as 
the rate of work due to the pressure. We emphasize that the 
continuous convective--diffusive equation displays very similar 
energy behavior (obviously the pressure term is absent there) 
\cite{EAk17}. This provides an additional indication of the 
similarity in the discrete setting.
\subsubsection*{Remark}
The transition from expression (\ref{eq:glb_flux}) to (\ref{eq:glb_flux alt}) 
is only possible in the continuous setting. In the discrete setting no 
closed-form expression for the Lagrange multiplier exists. This also 
applies to the localized version in Section \ref{sec:GE, subsec:LEE}.

\subsection{Local energy evolution}\label{sec:GE, subsec:LEE}
The procedure to find the local energy evolution is very similar to that 
of the global energy. Let $\omega \subset \Omega$ be an arbitrary 
subdomain with boundary $\pd \omega$, let $\Omega-\omega$ 
denote its complement and let their shared boundary denote 
$\chi_\omega = \pd \omega \cap \pd (\Omega-\omega)$. 
Figure \ref{fig:domain2} shows the subdomains and their boundaries.
\begin{figure}[h!]
  \begin{center}
    \begin{tikzpicture}
\draw[line width=0.3mm, black ] (2,0.4) .. controls (3.5,1) .. (5.1, -0.25) .. controls (6.5,-1.5) .. (5, -2.5) .. controls (4.5,-2.75) .. (3.5, -2).. controls (3.0,-1.5) .. (2.2, -2) .. controls (1.1,-2.5) .. (1.2, -1.0) .. controls (1.3,0.1) .. (2.0,0.4);
\draw[line width=0.3mm, black ] (5.1, -0.25) .. controls (3.6,-0.2) .. (3.5, -2);
\node[text width=3cm] at (3.5,-0.5) {$\Omega-\omega$};
\node[text width=3cm] at (6.0,-1.5) {$\omega$};
\node[text width=3cm] at (6.0,-3.0) {$\Gamma_h(\omega)$};
\node[text width=3cm] at (7.7,-1.0) {$\Gamma_g(\omega)$};
\node[text width=3cm] at (3.8,1.0) {$\Gamma_g(\Omega-\omega)$};
\node[text width=3cm] at (2.6,-2.8) {$\Gamma_h(\Omega-\omega)$};
\node[text width=3cm] at (5.3,-0.1) {$\chi_\omega$};
\node at (5cm,-1) {\pgfuseplotmark{cross*}};
\node[text width=3cm] at (7.0,-2.1) {$\boldsymbol{\backslash}$};
\node[text width=3cm] at (3.0,0.25) {$\boldsymbol{\backslash}$};
\node[text width=3cm] at (6.5, -0.25) {$\boldsymbol{/}$};
\node[text width=3cm] at (4.9, -2) {$\boldsymbol{/}$};
\end{tikzpicture}
  \end{center}
  \caption{Spatial domain $\Omega$ with a subdomain 
  $\omega\subset \Omega$. The shared boundary of $\omega$ 
  and its complement is $\chi_\omega$. The boundaries 
  $\Gamma_g$ and $\Gamma_h$ split according to $\omega$. 
  This is the same figure as in \cite{EAk17}.}
  \label{fig:domain2}
\end{figure}
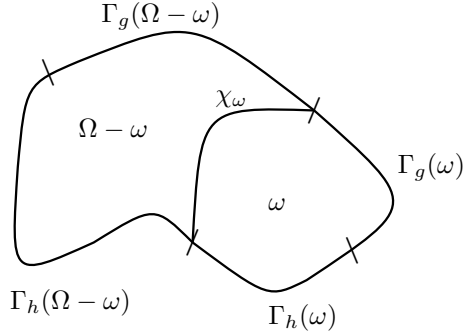\\

The continuity across the interface is enforced with a Lagrange multiplier 
in the appropriate space $\VV_\omega$. The discontinuous test function 
space writes $\WW_\omega$. The weak statement enforced on $\omega$ 
is again a mixed formulation and reads:\\

\textit{Find $\left(\left\{\mathbf{u}, p\right\}, \boldsymbol{\lambda}_{\omega}
\right) \in \WW \times \VV$ such that for all $\left(\left\{\mathbf{w}, q\right\}, \boldsymbol{\vartheta}\right) \in \WW \times \VV$,}
\begin{align}\label{sec:GE, subsec:LEE, standard weak form1}
\left(\mathbf{w},\boldsymbol{\lambda}_\omega\right)_{\chi_\omega} 
+\left(\mathbf{w},\boldsymbol{\lambda}_\omega\right)_{\Gamma_g \left(\omega\right)} 
&= B_{\omega,\Gamma_h}\left(\left\{\mathbf{u}, p\right\},\left\{\mathbf{w}, q\right\}\right)
- L_{\omega,\Gamma_h}(\left\{\mathbf{w},q\right\}), \nonumber \\
\left(\boldsymbol{\vartheta},[\![ \mathbf{u} ]\!] \right)_{\chi_\omega}
+\left(\boldsymbol{\vartheta},\mathbf{u}-\mathbf{g}\right)_{\Gamma_g\left(\omega\right)}&=0.
\end{align}
We have here employed the jump term $[\![ \mathbf{u} ]\!]$ given by
\begin{equation}
 [\![ \mathbf{u} ]\!]:=\mathbf{u}|_\omega -\mathbf{u}|_{\Omega-\omega},
\end{equation}
where the terms are defined on $\omega$ and $\Omega-\omega$, respectively. 
Furthermore, $\mathbf{n}_\omega$ is the outward normal of domain 
$\omega$, $u_{n_\omega}$ is the outward velocity in direction 
$\mathbf{n}_\omega$ and $\partial_{n_\omega}$ the direction derivative 
outward of $\omega$. The equivalence of this form with the strong form 
(\ref{sec:GE, subsec:SF, NS Strong}) leads to the expression of the Lagrange multiplier:
\begin{align}\label{}
  \boldsymbol{\lambda}_\omega= - u_{n_\omega} \mathbf{u} - p \mathbf{n}_\omega+\nu \pd_{n_\omega} \mathbf{u} ,
\end{align}
which is clearly the localized version of (\ref{sec:GE, subsec:SWF, LM}).
A direct consequence is the symmetry of the Lagrange multipliers 
(these are also called auxiliary fluxes in this setting, see \cite{HEML00}):
\begin{align}
  \boldsymbol{\lambda}_\omega+\boldsymbol{\lambda}_{\Omega-\omega}=\B{0},
\end{align}
i.e. that what flows out $\omega$ through $\chi_\omega$ enters its complement. 
The energy evolution linked to each of the domains is a natural split of the global energy evolution:
\begin{subequations}
  \label{gal: localized energy evolution}
  \begin{alignat}{1}
  \dfrac{{\rm d}}{{\rm d}t} E_\omega=&-\|\nu^{1/2}\nabla \mathbf{u}\|^2_{\omega}+ \left(\mathbf{u},\mathbf{f}\right)_{\omega}-\left(1,F_{\omega}\right)_{\partial \omega},
  \end{alignat}
\end{subequations}
with energy fluxes
\begin{align}
  \label{eq:glb_flux2}
F_{\omega} =
\left \{ \begin{array}{ll}
  -  \mathbf{g}\cdot \boldsymbol{\lambda}_\omega &\text{on}~~~ \Gamma_g\left(\omega\right),  \\
  |u_{n_\omega}|e  - \mathbf{u}\cdot\mathbf{h}&\text{on} ~~~ \Gamma_h\left(\omega\right),\\
    - \mathbf{u}\cdot \boldsymbol{\lambda}_\omega   &\text{on} ~~~\chi_\omega.
  \end{array} \right .
\end{align}
The last term  of (\ref{eq:glb_flux2}) redistributes energy over the domain. 
It represents an energy flux across the subdomain interface $\chi_\omega$ 
with a convective, a pressure and a viscous component. Similarly as before, 
substitution of the terms in the energy flux leads to
\begin{align}
  \label{eq:glb_flux alt2}
F_\omega = u_{n_\omega} (e+p) - \nu \pd_{n_\omega} e \quad \text{on} ~~~ \partial \omega.
\end{align}
This is obviously the localized version of (\ref{eq:glb_flux alt}).
\subsection*{Remark}
All statements of this Section are in the continuous setting. Hence, the 
standard discretization, i.e. the Galerkin method, displays the 
same correct energy behavior.

\subsection*{Remark}
The various boundary terms may distract the reader and do not 
contribute to the goal of this paper. Therefore we only consider 
boundary conditions precluding the energy flux $F$ on $\Gamma$. 
The homogeneous Dirichlet and periodic boundary conditions 
satisfy this purpose. Applying non-homogeneous boundaries 
is straightforward.\\

We continue this paper by discretizing the system according to 
the dynamic variational multiscale method with the target to 
closely resemble energy evolution 
(\ref{sec:GE, subsec:EE, energy1}) 
and (\ref{gal: localized energy evolution}). 

\section{Energy evolution of the variational multiscale method
 with dynamic small-scales}\label{sec:EESSM}
The convective--diffusive context \cite{EAk17} learns us that 
the dynamical structure of the small-scales is a requirement for 
the stabilized formulation to display the correct energy behavior. 
This allows to skip the static small-scales and to directly apply the 
dynamic modeling approach. We follow this road.
\subsection{The multiscale split}\label{sec:The multiscale split}
The variational multiscale split is nowadays a standard execution 
\cite{Hug95, Hug98} which we include here for the sake of 
completeness and notation. Employing the variational multiscale 
methodology the trial and weighting function spaces split into 
large- and small-scales as:
\begin{align}
  \WW = \WW^h \oplus \WW',
\end{align}
with $\WW^h$ and $\WW'$ containing the large-scales and
small-scales, respectively. The large-scale space is spanned 
by the finite dimensional numerical discretization while the 
fine-scales are its infinite dimensional complement. The fine-scale 
space $\WW'$ is also referred to as subgrid-scales since these 
scales are not reproduced by the grid. This decomposition 
implies the split of the solution and weighting functions as follows:
\begin{subequations}
  \begin{alignat}{2}
  \mathbf{U} =&~ \mathbf{U}^h + \mathbf{U}', \label{solution split}\\
  \mathbf{W} =&~ \mathbf{W}^h + \mathbf{W}',
  \end{alignat}
\end{subequations}
where $\mathbf{U}^h, \mathbf{W}^h \in \WW^h$ 
and $\mathbf{U}', \mathbf{W}' \in \WW'$ 
with $\mathbf{U}:=\left\{\mathbf{u},p\right\}, \mathbf{W}:=\left\{\mathbf{w},q\right\}$. 
Uniqueness follows when a projector 
$\mathscr{P}^h: \WW \rightarrow \WW^h$ is used for the splitting operation:
\begin{subequations}
  \label{eq:opt_proj}
  \begin{alignat}{1}
    \mathbf{U}^h &= \mathscr{P}^h \mathbf{U}, \label{scale sep proj}\\
    \mathbf{U}' &= \left(\mathscr{I}-\mathscr{P}^h\right) \mathbf{U},
  \end{alignat}
\end{subequations}
where $\mathscr{I}: \WW \rightarrow \WW$ is the identity operator. 
Employing both $\mathbf{W} = \mathbf{W}^h$ and 
$\mathbf{W} = \mathbf{W}'$, and the solution split 
(\ref{solution split}) in (\ref{sec:GE, subsec:SWF, standard weak form1})
 leads to the weak formulation:\\

\textit{Find $\mathbf{U}^h \in \WW^h,~\mathbf{U}' \in \WW'$ $\text{for all}~~\mathbf{W}^h \in \WW^h,~\mathbf{W}' \in \WW'$,}
\begin{subequations}
  \begin{alignat}{2}
  B_{\Omega}\left(\mathbf{U}^h+\mathbf{U}',\mathbf{W}^h\right) 
  =& \added[id=Rev.1]{L_{\Omega}}(\mathbf{W}^h)_{\Omega}, \quad \text{for all }\mathbf{W}^h\in \WW^h, \label{VMS large-scale eq} \\ 
  B_{\Omega}\left(\mathbf{U}^h+\mathbf{U}',\mathbf{W}'\right) 
  =& \added[id=Rev.1]{L_{\Omega}}(\mathbf{W}')_{\Omega}, \quad \text{for all }\mathbf{W}'\in \WW'. \label{VMS small-scale eq}
  \end{alignat}
\end{subequations}
Note that this is an infinite-dimensional system with unknowns 
$\mathbf{U}^h$ and $\mathbf{U}'$. Appropriately parameterizing 
the small-scales $\mathbf{U}'$ in terms of $\mathbf{U}^h$ 
converts (\ref{VMS large-scale eq}) into a solvable finite element problem. 
This conversion can be done with inspiration from (\ref{VMS small-scale eq}).
 For the technical details of the parameterization consult \cite{HugSan06}.
\subsection{Dynamic small-scales}\label{sec:Dynamic small-scales}
Here we employ the dynamic small-scales, see \cite{Cod02}, 
demanded by the convective--diffusive context for correct energy 
behavior \cite{EAk17}. The fine-scale model
\begin{align}\label{eq:Dynamic small-scales}
    \pd_t\left\{\hat{\mathbf{u}}', 0\right\}+\boldsymbol{\tau}^{-1}\left\{\hat{\mathbf{u}}', \hat{p}'\right\} 
    + \mathscr{R}\left(\left\{\mathbf{u}^h, p^h\right\}, \hat{\mathbf{u}}' \right)=0,
\end{align}
is an ordinary differential equation. The hat-sign is used to indicate 
a small-scale model instead of the actual small-scales. The intrinsic 
time scale $\boldsymbol{\tau}$ is a matrix of stabilization parameters, 
here $\boldsymbol{\tau} \in \mathbb{R}^{4\times 4}$,  with 
contributions for the two equations:
\begin{align}
    \boldsymbol{\tau}=\begin{pmatrix} \tau_M \B{I}_{3 \times 3} & \B{0}_3 \\ \B{0}_3^T & \tau_C \end{pmatrix}.
\end{align}
The local large-scale residual contains a momentum part 
$\B{r}_M$ and continuity part $r_C$ linked to the 
incompressibility constraint, respectively, given by
\begin{subequations}
  \label{vms: residuals}
  \begin{alignat}{1}
   \mathscr{R}\left(\left\{\mathbf{u}^h, p^h\right\}, \hat{\mathbf{u}}' \right)&=\left\{\B{r}_M(\left\{ \mathbf{u}^h, p^h\right\},\hat{\mathbf{u}}'), r_C(\mathbf{u}^h) \right\}^T,\\
   \B{r}_M &= \pd_t \mathbf{u}^h + \left(\left(\mathbf{u}^h+\hat{\mathbf{u}}'\right)\cdot\nabla\right) \mathbf{u}^h 
                  + \nabla p^h - \nu \Delta \mathbf{u}^h - \mathbf{f}, \label{vms: residuals mom}\\
   r_C &= \nabla \cdot \mathbf{u}^h. \label{vms: residuals cont}
  \end{alignat}
\end{subequations}
In the following we ignore the hat-sign. We employ a dynamic version 
of the stabilization parameters $\tau_M, \tau_C$ defined in \cite{BaCaCoHu07}. 
The details are provided in \ref{Appendix: Definition stabilization parameters}.
 The subscripts \textit{M} and \textit{C} refer to \textit{momentum} and \textit{continuity}, respectively. 
 Mirroring \cite{EAk17}, the momentum residual (\ref{vms: residuals mom}) uses 
 the full velocity $\mathbf{u}^h+\mathbf{u}'$. This creates a nonlinearity in the system. 
 Therefore we apply a standard iterative procedure to determine the small-scales.

Assume now that the domain $\Omega$ is partitioned into a set of elements $\Omega_e$. 
The domain of element interiors does not include the interior boundaries and denotes
\begin{equation}
\tilde{\Omega} = \displaystyle \bigcup_e \Omega_e.
\end{equation}

The resulting residual-based dynamic VMS weak formulation is\\

\textit{Find $\mathbf{U}^h \in \WW^h$ $\text{for all} ~~\mathbf{W}^h \in \WW^h$}
\begin{subequations}
  \label{standard dynamic VMS}
  \begin{alignat}{2}
     B_{\Omega}^{\text{VMSD}}\left(\mathbf{U}^h,\mathbf{W}^h\right)=& \added[id=Rev.1]{L_{\Omega}}(\mathbf{W}^h),
     \end{alignat}
\indent \textit{where}
\begin{alignat}{2}
     B_{\Omega}^{\text{VMSD}}\left(\mathbf{U}^h,\mathbf{W}^h\right) =& B_{\Omega}\left(\mathbf{U}^h,\mathbf{W}^h\right)
     + \left(\mathbf{w}^h, \pd_t \mathbf{u}'\right)_{\tilde{\Omega}}
     - \left( \nu \Delta \mathbf{w}^h, \mathbf{u}'\right)_{\tilde{\Omega}}  \nonumber \\
     &-\left(\nabla q^h, \mathbf{u}'\right)_{\tilde{\Omega}}-\left(\nabla \cdot \mathbf{w}^h, p'\right)_{\tilde{\Omega}}\nonumber \\
     & -\left(\nabla \mathbf{w}^h, \mathbf{u}^h\otimes \mathbf{u}'\right)_{\tilde{\Omega}} 
     -\left(\nabla \mathbf{w}^h, \mathbf{u}'\otimes \mathbf{u}^h\right)_{\tilde{\Omega}} 
     -\left(\nabla \mathbf{w}^h, \mathbf{u}'\otimes \mathbf{u}'\right)_{\tilde{\Omega}},\label{VMS bilinear}\\
&   \pd_t\left\{\mathbf{u}', 0\right\}+\boldsymbol{\tau}^{-1}\left\{\mathbf{u}', p'\right\} 
+ \mathscr{R}\left(\left\{\mathbf{u}^h, p^h\right\}, \mathbf{u}' \right)=0,\label{VMSD small scales}
  \end{alignat}
\end{subequations}
and where the additional D stands for \textit{dynamic}. When examining the last line of 
(\ref{VMS bilinear}), we recognize the following contributions. The first term is the 
SUPG contribution. The first two terms model the \textit{cross stress},  while the 
last term models the \textit{Reynolds stress}. Note that no spatial derivatives act 
on the small-scales. Furthermore, in contrast to static small-scales, the dynamic 
small-scale model (\ref{VMSD small scales}) is a separate equation and cannot 
directly be substituted into the large-scale equation (\ref{VMS bilinear}).

\subsection{Local energy evolution of the VMSD form}
\label{sec:Local energy evolution dynamic VMS}

To arrive at the local energy evolution of (\ref{standard dynamic VMS}), 
we extend the weak formulation to a Lagrange multiplier setting to allow 
discontinuous functions across subdomains, similar as 
(\ref{sec:GE, subsec:LEE, standard weak form1}). The weak statement, 
here stated for domain $\omega \subset \Omega$, reads \\

\textit{Find $\left(\mathbf{U}^h, \boldsymbol{\lambda}_{\omega}^h\right) \in \WW \times \VV$ 
such that for all $\left(\mathbf{W}^h, \boldsymbol{\vartheta}^h\right) \in \WW \times \VV$,}
\begin{subequations}
  \label{sec:VMS dyn LM form}
  \begin{alignat}{2}
\left(\mathbf{w}^h,\boldsymbol{\lambda}^h_\omega\right)_{\chi_\omega} 
&= B_{\omega}^{\text{VMSD}}\left(\mathbf{U}^h,\mathbf{W}^h\right)- \added[id=Rev.1]{L_{\omega}}(\mathbf{W}^h), \\
\left(\boldsymbol{\vartheta}^h,[\![\mathbf{u}^h]\!]\right)_{\chi_\omega}&=0,\\
  \pd_t\left\{\mathbf{u}', 0\right\}+\boldsymbol{\tau}^{-1}\left\{\mathbf{u}', p'\right\} 
  + \mathscr{R}\left(\left\{\mathbf{u}^h, p^h\right\}, \mathbf{u}' \right)&=0.\label{small scale dynamic VMS}
  \end{alignat}
\end{subequations}

To obtain the evolution of the \textit{local total energy} 
$E_\omega = \tfrac{1}{2}\left(\mathbf{u}^h + \mathbf{u}',\mathbf{u}^h + \mathbf{u}'\right)_{\tilde{\omega}}$ 
linked to the variational formulation (\ref{standard dynamic VMS}), we employ 
$\mathbf{w}^h = \mathbf{u}^h, q^h = p^h$ 
and $\boldsymbol{\vartheta}^h=\boldsymbol{\lambda}^h_{\omega} $ in (\ref{sec:VMS dyn LM form}). 
Adding $\mathbf{u}'$ times the momentum component of (\ref{small scale dynamic VMS}) 
integrated over $\tilde{\omega}$ eventually leads to
\begin{align}\label{energy evolution dyn VMS}
 \dfrac{{\rm d}}{{\rm d}t} E_\omega=& 
 - \|\nu^{1/2}\nabla \mathbf{u}^h\|_{\omega}^2
 + (\mathbf{u}^h,\mathbf{f})_{\omega}
 - (1,F_\omega^h)_{\chi_\omega}\nonumber  \\
                     & -\|\tau_M^{-1/2}\mathbf{u}'\|_{\tilde{\omega}}^2 
                     + (\mathbf{u}',\mathbf{f})_{\tilde{\omega}} 
                     + 2(\nu \Delta \mathbf{u}^h,\mathbf{u}')_{\tilde{\omega}}\nonumber \\
                     &+(\nabla \cdot \mathbf{u}^h,p')_{\tilde{\omega}}
                     +\left(\nabla \mathbf{u}^h , (\mathbf{u}^h 
                     + \mathbf{u}') \otimes (\mathbf{u}^h 
                     + \mathbf{u}') \right)_{\tilde{\omega}} 
                     - \left(\mathbf{u}',(\mathbf{u}^h
                     +\mathbf{u}')\cdot \nabla \mathbf{u}^h \right)_{\tilde{\omega}},
\end{align}
where
\begin{equation}\label{energy flux VMS}
  F_\omega^h= - \boldsymbol{\lambda}^h_\omega \cdot \mathbf{u}^h.
\end{equation}
The first line closely resembles the continuous energy evolution relation. 
Each one of the other terms appears as a result of the VMS stabilization. 
The first term of the second line represents the numerical dissipation due 
to the missing small-scales. This contributes to a decay of the energy, 
which is favorable from a stability argument. The second term is the 
power exerted by the body force on the small-scales, this term closely 
resembles its large-scale counterpart. The remaining terms have no 
continuous counterpart. With the current small-scale model, the 
small-scale pressure term dissipates energy\footnote{The small-scale pressure expression can be substituted into this 
term to arrive at $(\nabla \cdot \mathbf{u}^h,p')_{\tilde{\omega}} = -|| \tau_C^{-1/2} p'||^2_{\tilde{\omega}}$.
Note that it vanishes when employing a divergence-conforming discrete velocity space.}. 
The signs of the other terms are undetermined and therefore these can create energy artificially. 
The term $2(\nu \Delta \mathbf{u}^h,\mathbf{u}')_{\tilde{\omega}}$ can be bounded by both the 
physical dissipation $ \|\nu^{1/2}\nabla \mathbf{u}^h\|_{\omega}^2$ and numerical dissipation 
$\|\tau_M^{-1/2}\mathbf{u}'\|_{\tilde{\omega}}^2$ using a standard argument. However, this 
results in an overall dissipation that can be smaller than the physical one. This is deemed 
undesirable. Note that it is comparable with that of the dynamic VMS stabilized form in the 
convective--diffusive context. The contrast occurs in the last line which is linked to the 
incompressibility constraint (\ref{sec:GE, subsec:SF, NS Strong, div free eq}) and the 
small-scale pressure. Inspired by the convective--diffusive context, the next Section 
rectifies the method to closely resemble the energy behavior of the continuous setting.

\subsection*{Remark}
Employing $\omega=\Omega$, and hence $\tilde{\omega}=\tilde{\Omega}$, 
provides the global energy evolution of (\ref{standard dynamic VMS}):
\begin{align}\label{energy evolution dyn VMS global}
 \dfrac{{\rm d}}{{\rm d}t} E_\Omega=
 & - \|\nu^{1/2}\nabla \mathbf{u}^h\|_{\Omega}^2+ (\mathbf{u}^h,\mathbf{f})_{\Omega}\nonumber  \\
 & -\|\tau_M^{-1/2}\mathbf{u}'\|_{\tilde{\Omega}}^2 + (\mathbf{u}',\mathbf{f})_{\tilde{\Omega}} 
 + 2(\nu \Delta \mathbf{u}^h,\mathbf{u}')_{\tilde{\Omega}}\nonumber \\
  &+(\nabla \cdot \mathbf{u}^h,p')_{\tilde{\Omega}}
    +\left(\nabla \mathbf{u}^h , (\mathbf{u}^h + \mathbf{u}') \otimes (\mathbf{u}^h + \mathbf{u}') \right)_{\tilde{\Omega}} 
    - \left(\mathbf{u}',(\mathbf{u}^h+\mathbf{u}')\cdot \nabla \mathbf{u}^h \right)_{\tilde{\Omega}}.
\end{align}
\section{Toward a stabilized formulation with correct energy behavior}
\label{sec:Towards correct energy behavior}
This Section presents the procedure to remedy the incorrect energy behavior 
(\ref{energy evolution dyn VMS}) of the dynamic VMS formulation (\ref{standard dynamic VMS}). 
The first ingredient is the switch from the conservative form to a skew-symmetric form 
with the help of the divergence-free velocity field constraint. Next, we employ the natural 
choice of a Stokes-projector and demand divergence-free small-scales. In view of the 
convective--diffusive context, we use $H_0^1$ small-scales to treat the small-scale viscous term.

\subsection{Design condition}\label{subsec: Design condition}
We present a design condition which clarifies the desirable energy behavior of the formulation. 
The variational weak formulation corresponding to (\ref{sec:GE, subsec:SF, NS Strong}) 
is demanded to satisfy the local energy behavior:
\begin{align}\label{sec:EESSM, subsec: EE, energy evo design cond}
  \dfrac{{\rm d}}{{\rm d}t} E_{\omega}= & - \|\nu^{1/2}  \nabla \mathbf{u}^h \|_{\omega}^2 
  + (\mathbf{u}^h ,\mathbf{f})_{\omega}-(1,F_\omega^h)_{\chi_\omega} \nonumber\\
                    & - \|\tau_M^{-1/2} \mathbf{u}'\|_{\tilde{\omega}}^2 + (\mathbf{u}' ,\mathbf{f})_{\tilde{\omega}},
\end{align}
with \textit{exact} divergence-free velocity fields. Note that this requirement is very 
similar to that of the convective--diffusive context \cite{EAk17} where the convective velocity is assumed solenoidal.

\added[id=Rev.2]{
\subsection*{Remark}
In the following we use the ingredients mentioned above to convert the VMS formulation (\ref{sec:VMS dyn LM form}) into a method that satisfies the design condition.
It is important to realize that the small-scales employed in the formulation are determined by a model equation.
This implies that these properties are not necessarily valid for the \textit{model} small-scales.
In contrast, the \textit{exact} small-scales do satisfy these properties.
The model small-scales approximate its exact counterpart which justifies the judicious use of these properties to construct a method that satisfies the design condition.
}
\subsection{Skew-symmetric form}
\label{subsec:CF}

We employ a multiscale form of the skew-symmetric formulation 
(see (\ref{sec:GE, subsec:SWF, standard weak form1: convective})) 
to eliminate the convective contributions in (\ref{energy evolution dyn VMS}). 
Considering the convective terms in isolation, we cast them into the following form:

\begin{align}\label{convective derivation}
 -(\nabla \mathbf{w}^h, (\mathbf{u}^h+\mathbf{u}')\otimes (\mathbf{u}^h+\mathbf{u}'))_{\tilde{\Omega}}
=&- \left( \left(\mathbf{u}^h  + \mathbf{u}'\right)\cdot \nabla \mathbf{w}^h  ,  \mathbf{u}^h \right)_{\tilde{\Omega}}
- \left( \left(\mathbf{u}^h  + \mathbf{u}'\right)\cdot \nabla \mathbf{w}^h  ,  \mathbf{u}' \right)_{\tilde{\Omega}} \nonumber\\
 =& \tfrac{1}{2} \left(\mathbf{w}^h, \left(\mathbf{u}^h  + \mathbf{u}'\right)\cdot \nabla \mathbf{u}^h\right)_{\tilde{\Omega}}
- \tfrac{1}{2} \left(\left(\mathbf{u}^h  + \mathbf{u}'\right)\cdot \nabla \mathbf{w}^h, \mathbf{u}^h \right)_{\tilde{\Omega}} \nonumber \\
 &+\tfrac{1}{2}\left(\mathbf{u}^h,\mathbf{w}^h\nabla \cdot \left(\mathbf{u}^h+\mathbf{u}'\right)\right)_{\tilde{\Omega}}
 - \left( \left(\mathbf{u}^h  + \mathbf{u}'\right)\cdot \nabla \mathbf{w}^h  ,  \mathbf{u}' \right)_{\tilde{\Omega}} \nonumber\\
  =& \tfrac{1}{2} \left(\mathbf{w}^h, \left(\mathbf{u}^h  + \mathbf{u}'\right)\cdot \nabla \mathbf{u}^h\right)_{\tilde{\Omega}}
- \tfrac{1}{2} \left(\left(\mathbf{u}^h  + \mathbf{u}'\right)\cdot \nabla \mathbf{w}^h, \mathbf{u}^h \right)_{\tilde{\Omega}} \nonumber \\
 &- \left( \left(\mathbf{u}^h  + \mathbf{u}'\right)\cdot \nabla \mathbf{w}^h  ,  \mathbf{u}' \right)_{\tilde{\Omega}},
\end{align}
where we have employed the multiscale incompressibility constraint 
$\nabla \cdot \mathbf{u}=\nabla \cdot (\mathbf{u}^h+\mathbf{u}')=0$ 
in the last equality. The last expression is incorporated into the formulation. 
The resulting residual-based skew-symmetric VMS weak formulation is\\

\textit{Find $\mathbf{U}^h \in \WW^h$ $\text{such that for all} ~~\mathbf{W}^h \in \WW^h$,}
\begin{subequations}
  \label{sec:Towards correct energy, VMSC}
  \begin{alignat}{2}
     C_{\Omega}^{\text{VMSD}}\left(\mathbf{U}^h,\mathbf{W}^h\right)=& \added[id=Rev.1]{L_{\Omega}}(\mathbf{W}^h),
     \end{alignat}
\indent \textit{where}
\begin{alignat}{2}
     C_{\Omega}^{\text{VMSD}}\left(\mathbf{U}^h,\mathbf{W}^h\right) =& 
     C_{\Omega}\left(\mathbf{U}^h,\mathbf{W}^h\right)+ \left(\mathbf{w}^h, \pd_t \mathbf{u}'\right)_{\tilde{\Omega}}
     - \left( \nu \Delta \mathbf{w}^h, \mathbf{u}'\right)_{\tilde{\Omega}}\nonumber\\
     &-\left(\nabla q^h, \mathbf{u}'\right)_{\tilde{\Omega}}-\left(\nabla \cdot \mathbf{w}^h, p'\right)_{\tilde{\Omega}}\nonumber \\
     & +\tfrac{1}{2} \left(\mathbf{w}^h, \mathbf{u}' \cdot \nabla \mathbf{u}^h\right)_{\tilde{\Omega}}
- \tfrac{1}{2} \left( \mathbf{u}'\cdot \nabla \mathbf{w}^h, \mathbf{u}^h \right)_{\tilde{\Omega}} \nonumber \\
 &- \left( \left(\mathbf{u}^h  + \mathbf{u}'\right)\cdot \nabla \mathbf{w}^h  ,  \mathbf{u}' \right)_{\tilde{\Omega}},
 \label{sec:Towards correct energy, VMSC, bilinear form}\\
 \pd_t\left\{\mathbf{u}', 0\right\}+\boldsymbol{\tau}^{-1}\left\{\mathbf{u}', p'\right\} + \mathscr{R}\left(\left\{\mathbf{u}^h, p^h\right\}, \mathbf{u}' \right)&=0.
  \end{alignat}
\end{subequations}
This eliminates the convective contributions from the local energy evolution equation:
\begin{align}\label{sec:Towards correct energy, VMSC: energy}
 \dfrac{{\rm d}}{{\rm d}t} E_\omega=& 
 - \|\nu^{1/2}\nabla \mathbf{u}^h\|_{\omega}^2
 + (\mathbf{u}^h,\mathbf{f})_{\omega}- (1,F_\omega^h)_{\chi_\omega}\nonumber  \\
& -\|\tau_M^{-1/2}\mathbf{u}'\|_{\tilde{\omega}}^2 
+ (\mathbf{u}',\mathbf{f})_{\tilde{\omega}} + 2(\nu \Delta \mathbf{u}^h,\mathbf{u}')_{\tilde{\omega}} 
+(\nabla \cdot \mathbf{u}^h,p')_{\tilde{\omega}}.
\end{align}


\subsection{Stokes projector}
In the convective--diffusive context a $H_0^1$-orthogonality of the small-scale viscous term 
is required for correct energy behavior. This is the distinguished limit of $Pe \rightarrow 0$ of 
the steady convection--diffusion equations, where $Pe$ is the P\'{e}clet number. 
Its Navier--Stokes counterpart is to apply a Stokes-projector which is based on the distinguished 
limit $Re \rightarrow 0$ of the steady incompressible Navier--Stokes equations. 
Here $Re$ is the Reynolds number. Thus, applying a Stokes projection on the large-scale 
equation seems a natural choice. Moreover, it reduces the variational form in the limit 
$Re \rightarrow 0$ to the standard Galerkin method. This is a valid and established method 
in that regime, provided compatible discretizations for the velocity and pressure spaces are used. 

For the scale separation (\ref{eq:opt_proj}) we select the Stokes projector given by\\

$\mathscr{P}_{\text{Stokes}}^h: \mathbf{U} \in \WW \rightarrow \mathbf{U}^h \in \WW^h$: \textit{Find $\mathbf{U}^h \in \mathcal{W}^h$ 
such that for all $\mathbf{W}^h \in \WW^h$,}
\begin{subequations}\label{Stokes projector}
  \begin{alignat}{2}
    \left(\nu \Delta \mathbf{w}^h, \mathbf{u}^h\right)_\Omega + \left(\nabla \cdot \mathbf{w}^h,p^h\right)_\Omega
    &= \left(\nu \Delta \mathbf{w}^h, \mathbf{u}\right)_\Omega + \left(\nabla \cdot \mathbf{w}^h,p\right)_\Omega,\\
    \left(\nabla q^h, \mathbf{u}^h \right)_\Omega &= \left(\nabla q^h, \mathbf{u} \right)_\Omega,
  \end{alignat}
\end{subequations}
in the bilinear form (\ref{sec:Towards correct energy, VMSC, bilinear form}). \added[id=Rev.2]{Note that this projector only makes sense if the elements of $\mathcal{W}^h$ are inf--sup stable and the velocities are at least $C^1$-continuous. The numerical results presented in Section \ref{sec:ns_case} fulfill this requirement: quadratic NURBS basis functions are employed. However, note that the final form, given in \ref{Appendix: Galerkin/least-squares formulation with dynamic divergence-free small-scales}, does not have the smoothness restriction.}

As a consequence we assume the modeled small-scales to satisfy the 
orthogonality induced by the Stokes operator:
\begin{subequations}\label{Stokes projector orthog}
  \begin{alignat}{2}
    \left(\nu \Delta \mathbf{w}^h, \mathbf{u}'\right)_{\tilde{\Omega}} + \left(\nabla \cdot \mathbf{w}^h,p'\right)_{\tilde{\Omega}}&= 0,\\
    \left(\nabla q^h, \mathbf{u}' \right)_{\tilde{\Omega}} &= 0,
  \end{alignat}
\end{subequations}
for all $\mathbf{W}^h \in \WW^h$ . This converts (\ref{sec:Towards correct energy, VMSC}) 
into the simplified formulation:\\

\textit{Find $\mathbf{U}^h \in \WW^h$ $\text{such that for all} ~~\mathbf{W}^h \in \WW^h$}
\begin{subequations}
  \label{VMSDS form}
  \begin{alignat}{2}
     S_{\Omega}\left(\mathbf{U}^h,\mathbf{W}^h\right)=& \added[id=Rev.1]{L_{\Omega}}(\mathbf{W}^h),\label{large scale eq VMSDS} 
     \end{alignat}
\indent \textit{where}
\begin{alignat}{2}
     S_{\Omega}\left(\mathbf{U}^h,\mathbf{W}^h\right) =& 
     C_{\Omega}\left(\mathbf{U}^h,\mathbf{W}^h\right)
     + \left(\mathbf{w}^h, \pd_t \mathbf{u}'\right)_{\tilde{\Omega}}\nonumber \\
     & +\tfrac{1}{2} \left(\mathbf{w}^h, \mathbf{u}' \cdot \nabla \mathbf{u}^h\right)_{\tilde{\Omega}}
- \tfrac{1}{2} \left( \mathbf{u}'\cdot \nabla \mathbf{w}^h, \mathbf{u}^h \right)_{\tilde{\Omega}} \nonumber \\
 &- \left( \left(\mathbf{u}^h  + \mathbf{u}'\right)\cdot \nabla \mathbf{w}^h  ,  \mathbf{u}' \right)_{\tilde{\Omega}},\\
   \pd_t \mathbf{u}'+\tau_M^{-1} \mathbf{u}'  + \B{r}_M =& 0, \label{small scale eq VMSDS}
  \end{alignat}
\end{subequations}
where the $S$ abbreviates \textit{Stokes}. Note that the small-scale pressure terms 
have vanished from the formulation. The energy linked to this formulation is
\begin{align}\label{energy evolution VMSDS form}
 \dfrac{{\rm d}}{{\rm d}t} E_\omega=& 
 - \|\nu^{1/2}\nabla \mathbf{u}^h\|_{\omega}^2
 + (\mathbf{u}^h,\mathbf{f})_{\omega}- (1,F_\omega^h)_{\chi_\omega}\nonumber  \\
& -\|\tau_M^{-1/2}\mathbf{u}'\|_{\tilde{\omega}}^2 + (\mathbf{u}',\mathbf{f})_{\tilde{\omega}} 
+ (\nu \Delta \mathbf{u}^h,\mathbf{u}')_{\tilde{\omega}} -(\nabla p^h,\mathbf{u}')_{\tilde{\omega}}.
\end{align}
To fulfill the design condition (\ref{sec:EESSM, subsec: EE, energy evo design cond}), 
the last two terms of (\ref{energy evolution VMSDS form}) need to be eliminated, i.e.
\begin{align}\label{eq: necessary optimality projector}
 (\nu \Delta \mathbf{u}^h,\mathbf{u}')_{\tilde{\Omega}} -(\nabla p^h,\mathbf{u}')_{\tilde{\Omega}} = 0.
\end{align}
 There are various options available  to accomplish this. 
Before sketching some of these options we first like to note the following.
Augmenting the undesirable terms of (\ref{energy evolution VMSDS form}) 
with $(\nabla \cdot \mathbf{u}^h,p')$ results in the requirement
\begin{align}
 (\nu \Delta \mathbf{u}^h,\mathbf{u}')_{\tilde{\omega}} - (\nabla p^h,\mathbf{u}')_{\tilde{\omega}}+(\nabla \cdot \mathbf{u}^h,p')_{\tilde{\omega}} =0.
\end{align}
This is a well-defined orthogonality induced by the Stokes operator, given in (\ref{Stokes projector orthog}).
The augmented term would appear if $\nabla p'$ in the small-scale momentum 
equation is not neglected\footnote{
Including the small-scale pressure in the residual augments the right-hand 
side of (\ref{energy evolution VMSDS form}) with the term $\left(\nabla p',\mathbf{u}'\right)$. 
Next, by using the strong form continuity equation weighted with the 
small-scale pressure, i.e. $\left(p',\nabla \cdot (\mathbf{u}^h+\mathbf{u}')\right)=0$, 
this term converts into $(\nabla \cdot \mathbf{u}^h,p')$.
}.
Note that this term is not (easily) computable and therefore usually omitted in the formulation.

The required orthogonality (\ref{eq: necessary optimality projector}) can be either 
\textit{assumed} or \textit{enforced} \cite{EAk17}. We discuss four options here.
\begin{itemize}
\item First we could assume the orthogonality in the small-scale equation (\ref{small scale eq VMSDS}).
This orthogonality has previously been assumed to modify the large-scale equation (\ref{large scale eq VMSDS}). 
Assuming it in the small-scale equation results in a stable method with the desired energy property. 
However the small-scale model is not residual-based anymore. This results in an inconsistent method. 
We do not further consider this option.
\item Alternatively, we could assume the orthogonality in the large-scale equation 
(\ref{large scale eq VMSDS}) again. This converts the formulation into a GLS method. 
This method includes a PSPG term, $-(\nabla q^h,\mathbf{u}')_{\tilde{\Omega}}$, and therefore 
pointwise divergence-free solutions cannot be guaranteed. The formulation harms the design 
condition of Section \ref{subsec: Design condition} and is therefore omitted.
\item Another option is to enforce the required orthogonality using Lagrange-multipliers. 
This is not straightforward and is deemed unnecessarily expensive.
\item The path we propose is to cure the unwanted terms separately by combining 
the second and third options. The approach is to (i) enforce divergence-free small-scales 
to eliminate the second term of \added[id=Rev.2]{(\ref{eq: necessary optimality projector})} and (ii) assume an 
$H_0^1$-orthogonality to erase the first term of \added[id=Rev.2]{(\ref{eq: necessary optimality projector})}. 
Sections \ref{subsec:Divergence-free small-scales} and \ref{subsec: H_0^1-orthogonal small-scales} 
respectively describe these steps. 
\end{itemize}

\subsection{Divergence-free small-scales}\label{subsec:Divergence-free small-scales}
The last term of (\ref{energy evolution VMSDS form}) disappears when enforcing 
divergence-free small-scales. We handle this with a projection operator on the small-scales:\\

$\mathscr{P}^h_{\text{div}}: \mathbf{U} \in \WW \rightarrow \mathbf{U}^h \in \WW^h$: 
\textit{Find $\mathbf{U}^h \in \mathcal{W}^h$ such that for all $\mathbf{W}^h \in \WW^h$,}
\begin{align}\label{div projection}
 \left(\nabla q^h,\mathbf{u}^h \right)_{\Omega}&= \left(\nabla q^h,\mathbf{u} \right)_{\Omega},
\end{align}
with corresponding orthogonality:
\begin{align}\label{div projection}
 \left(\nabla q^h,\mathbf{u}' \right)_{\tilde{\Omega}}&= 0, \quad \text{for all}\quad  \mathbf{W}^h \in \WW^h.
\end{align}
This orthogonality defines the fine-scale space $\WW'$ which represents the 
orthogonal component of $\WW^h$ in terms of the projection (\ref{div projection}) as
\begin{equation}\label{sec:EEDF, subsec:OSS, constricted fine-scale space}
  \begin{array}{l r l}
    \WW'=\WW'_{\text{div}}:=\Big\{\left\{\mathbf{u},p\right\} \in \WW;  
    & \left. \left(\nabla \theta^h,\mathbf{u} \right)_{\Omega}=0, \right. 
    &  \text{for all } \theta^h \in \mathcal{P}^h \Big\},
  \end{array}
\end{equation}
where the space $\mathcal{P}^h$ is \added[id=Rev.1]{the} pressure part of $\WW^h=\mathcal{U}^h\times \mathcal{P}^h$. 
Directly employing this divergence-free space indeed eliminates the last term of (\ref{energy evolution VMSDS form}). 
However the small-scale solution space is infinite dimensional, and therefore not applicable in the numerical method. 
As before, we avoid dealing with this space by using a Lagrange-multiplier construction yielding a mixed formulation.
 Opening the solution space leads to the formulation:\\

\textit{Find $\left(\mathbf{U}^h, \zeta^h\right) \in \WW^h\times \mathcal{P}^h$ such that for all $\left(\mathbf{W}^h,\theta^h \right) \in \WW^h\times\mathcal{P}^h$,}
\begin{subequations}
  \label{sec:Towards correct energy, GLS}
  \begin{alignat}{2}
     S_{\Omega}^{\text{div}}\left(\left(\mathbf{U}^h, \zeta^h\right),\left(\mathbf{W}^h, \theta^h\right) \right) =&
 \added[id=Rev.1]{L_{\Omega}}(\mathbf{W}^h)_{\Omega}, 
     \label{sec:Towards correct energy, GLS, large scale}
     \end{alignat}
\indent \textit{where}
\begin{alignat}{2}
     S_{\Omega}^{\text{div}}\left(\left(\mathbf{U}^h, \zeta^h\right),\left(\mathbf{W}^h, \theta^h\right) \right)=&
      S_{\Omega}\left(\mathbf{U}^h,\mathbf{W}^h\right)+\left( \nabla \theta^h, \mathbf{u}'\right)_{\tilde{\Omega}},\\
  \pd_t \mathbf{u}'+\tau_M^{-1} \mathbf{u}' + \nabla \zeta^h  + \B{r}_M =& 0.
  \end{alignat}
\end{subequations}
Obviously, this form follows the energy evolution
\begin{align}\label{energy evolution VMSDSO form}
 \dfrac{{\rm d}}{{\rm d}t} E_\omega=& - \|\nu^{1/2}\nabla \mathbf{u}^h\|_{\omega}^2
 + (\mathbf{u}^h,\mathbf{f})_{\omega}- (1,F_\omega^h)_{\chi_\omega}\nonumber  \\
 & -\|\tau_M^{-1/2}\mathbf{u}'\|_{\tilde{\omega}}^2 
 + (\mathbf{u}',\mathbf{f})_{\tilde{\omega}} 
 + (\nu \Delta \mathbf{u}^h,\mathbf{u}')_{\tilde{\omega}}.
\end{align}
\added[id=Both Rev]{
\subsection*{Remark}
Note that enforcing divergence-free small-scales has introduced an additional equation in the system. The new method has $5$ global variables instead of $4$ leading to a commensurate increase in computational time. The added block diagonal term is a diffusion matrix which does not further complicate the saddle point structure of the problem.}

\subsection{$H_0^1$-orthogonal small-scales}
\label{subsec: H_0^1-orthogonal small-scales}
In the energy evolution (\ref{energy evolution VMSDSO form}) unwanted artificial energy 
can only be created by the term $\left(\nu \Delta \mathbf{u}^h,\mathbf{u}'\right)_{\tilde{\omega}}$. 
Employing the orthogonality induced by the $H_0^1$-seminorm,
\begin{align}
(\nu \Delta \mathbf{w}^h,\mathbf{u}')_{\tilde{\Omega}} = 0 \quad \text{for~all}~\mathbf{W}^h \in \WW^h,
\end{align}
obviously cancels this term. To avoid dealing with a larger system of equations, we do not enforce the 
orthogonality but we assume it in the large-scale equation (\ref{sec:Towards correct energy, GLS, large scale}). 
This leads to a consistent GLS method. The resulting GLSDD-formulation reads:\\

\textit{Find $\left(\mathbf{U}^h, \zeta^h\right) \in \WW^h\times \mathcal{P}^h$ such that for all $\left(\mathbf{W}^h,\theta^h \right) \in \WW^h\times\mathcal{P}^h$,}
\begin{subequations}
  \label{sec:Towards correct energy, GLSDSD}
  \begin{alignat}{2}
     S_{\Omega}^{\text{GLSDD}}\left(\left(\mathbf{U}^h, \zeta^h\right),\left(\mathbf{W}^h, \theta^h\right) \right)=& \added[id=Rev.1]{L_{\Omega}}(\mathbf{W}^h),
     \end{alignat}
\indent \textit{where}
\begin{alignat}{2}
     S_{\Omega}^{\text{GLSDD}}\left(\left(\mathbf{U}^h, \zeta^h\right),\left(\mathbf{W}^h, \theta^h\right) \right) =&
      S_{\Omega}^{\text{div}}\left(\left(\mathbf{U}^h, \zeta^h\right),\left(\mathbf{W}^h, \theta^h\right)\right)
      + \left( \nu \Delta \mathbf{w}^h, \mathbf{u}'\right)_{\tilde{\Omega}},\nonumber \\
  \pd_t \mathbf{u}'+\tau_M^{-1} \mathbf{u}' + \nabla \zeta^h  + \B{r}_M =& 0.
  \end{alignat}
\end{subequations}
In the abbreviation GLSDD we follow the same structure as before where the last two D's 
stand for \textit{dynamic, divergence-free small-scale velocities}\footnote{\added[id=Rev.2]{The name GLS refers to the convection--diffusion part of the problem.}}. 
This method displays the correct-energy behavior:
\begin{equation}\label{energy evol GLSDD}
\begin{array}{l l}
 \dfrac{{\rm d}}{{\rm d}t} E_{\omega}= & - \|\nu^{1/2}  \nabla \mathbf{u}^h \|_{\omega}^2 + (\mathbf{u}^h ,\mathbf{f})_{\omega}-(1,F_\omega^h)_{\chi_\omega}\\
                    & - \|\tau_M^{-1/2} \mathbf{u}'\|_{\tilde{\omega}}^2 + (\mathbf{u}' ,\mathbf{f})_{\tilde{\omega}}.
\end{array}
\end{equation}
The full expansion of this novel formulation is included in 
\ref{Appendix: Galerkin/least-squares formulation with dynamic divergence-free small-scales} for clarity. 

\subsection{Local energy backscatter}
The separate energy evolution of the large- and small-scales deduces in a similar fashion as above. 
The large-scale energy $E^h_\omega=\tfrac{1}{2}(\mathbf{u}^h,\mathbf{u}^h)_{\omega}$ and the 
small-scale energy $E'_\omega=\tfrac{1}{2}(\mathbf{u}',\mathbf{u}')_{\tilde{\omega}}$ do not add 
up to the total energy $E_{\omega}$ because of the missing cross terms. This energy is stored 
in \added[id=Rev.1]{an} intermediate (buffer) regime which we denote with 
$E^{h'}_\omega=(\mathbf{u}^h,\mathbf{u}')_{\tilde{\omega}}$. The energy evolution takes the form:
\begin{subequations}\label{eq: local energy backscatter}
\begin{alignat}{2}
  \dfrac{{\rm d}}{{\rm d}t} E^h_\omega  & =
   - \| \nu^{1/2} \nabla \mathbf{u}^h \|^2_{\omega} 
   + \left(\mathbf{u}^h,\mathbf{f}\right)_{\omega} -(1,F_\omega)_{\chi_\omega} 
   + \left(\left(\mathbf{u}^h+\mathbf{u}'\right)\cdot \nabla \mathbf{u}^h,\mathbf{u}'\right)_{\tilde{\omega}} -( \mathbf{u}^h, \partial_t  \mathbf{u}' )_{\tilde{\omega}},  \\
  \dfrac{{\rm d}}{{\rm d}t} E^{h'}_\omega  & = \left(\mathbf{u}^h,\partial_t \mathbf{u}'\right)_{\tilde{\omega}} 
  + \left(\mathbf{u}',\partial_t \mathbf{u}^h\right)_{\tilde{\omega}},   \\
  \dfrac{{\rm d}}{{\rm d}t} E'_\omega   & = - \| \tau_M^{-1/2} \mathbf{u}'\|^2_{\tilde{\omega}} 
  + \left(\mathbf{u}' ,\mathbf{f}\right)_{\tilde{\omega}} 
  - \left(\left(\mathbf{u}^h+\mathbf{u}'\right)\cdot \nabla \mathbf{u}^h,\mathbf{u}'\right)_{\tilde{\omega}} 
  -( \mathbf{u}',  \partial_t \mathbf{u}^h )_{\tilde{\omega}}.
  \end{alignat}
\end{subequations}
The result mirrors to the convective--diffusive context with as convective velocity now the 
total velocity $\mathbf{u}^h+\mathbf{u}'$. There is a direct exchange of convective 
energy between the large-scale and small-scales. Clearly the superposition of 
(\ref{eq: local energy backscatter}) yields (\ref{energy evol GLSDD}).
\subsection{Time-discrete energy behavior}\label{subsec:Time-discrete energy behavior}
The generalized-$\alpha$ method serves as time-integrator. Mirroring the convective--diffusive 
context \cite{EAk17}, and using the same notation, we eventually obtain for $\alpha_m = \gamma$:
\begin{align}\label{discretized energy}
  E_{n+1}=E_n-\Delta t^2 (\alpha_f-\onehalf) \| \dot{\mathbf{u}}_{n+\alpha_m} \|^2_{\Omega}&
  -  \Delta t\| \nu^{1/2}\nabla \mathbf{u}^h_{n+\alpha_f} \|^2_{\Omega}
  -\Delta t\| \tau_{\text{dyn}}^{-1/2} \mathbf{u}_{n+\alpha_f}'\|^2_{\tilde{\Omega}} \nonumber \\
  &+ \Delta t(\mathbf{u}^h_{n+\alpha_f} ,f)_{\Omega}+ \Delta t(\mathbf{u}_{n+\alpha_f}' ,f)_{\tilde{\Omega}}.
\end{align}
Hence, we have a decay of the discretized energy when, in absence of forcing, 
$\alpha_f \geq \tfrac{1}{2}$. In the numerical implementation we use 
$\alpha_f = \alpha_m = \gamma = \tfrac{1}{2}$ for the stability and 
second-order accuracy properties \cite{Hul93}.

\section{Conservation properties}
\label{sec:CP}

Conservation of physical quantities in the numerical formulation is an often 
sought-after property. In this Section we derive the various conservation 
properties (continuity, linear momentum, angular momentum) of the proposed 
formulation (\ref{sec:Towards correct energy, GLSDSD}). We \added[id=Rev.1]{prove} these by 
selecting the appropriate weighting functions. The conservation properties 
hold on both a global and a local scale. Therefore we omit the domain 
subscript in the following.

\subsection{Continuity}
\label{sec:CP, subsec:C}

Employing the weighting function $\mathbf{w}^h=\mathbf{0},~\theta^h=0$ 
in (\ref{sec:Towards correct energy, GLSDSD}) yields
\begin{align}
 (q^h,\nabla \cdot \mathbf{u}^h) =  0.
\end{align}
The choice $q^h= \nabla\cdot \mathbf{u}^h$ \added[id=Rev.1]{proves} the \textit{pointwise} satisfaction of incompressibility constraint\footnote{\added[id=Rev.1]{Note that in general this weighting function choice is not allowed. We employ the IGA spaces with stable velocity and pressure pairs that do allow this choice.}} 
\begin{align}\label{div free}
 ||\nabla \cdot \mathbf{u}^h||^2=0 \quad \Rightarrow \quad \nabla \cdot \mathbf{u}^h=0 \quad \text{for all} \quad \mathbf{x} \in \Omega.
\end{align}
\added[id=Authors]{Furthermore,} the choice of weighting functions $\mathbf{w}^h=\mathbf{0},~q^h=0$ leads to divergence-free small-scale velocities in the following sense:
\begin{align}
 (\nabla \theta^h,  \mathbf{u}') =  0.
\end{align}

\subsection{Linear momentum}
\label{sec:CP, subsec:GLM}

We substitute the weighting functions 
$\left(\mathbf{w}^h,q^h, \theta^h\right)=\left(\mathbf{e}_i,0, -\tfrac{1}{2}\mathbf{e}_i\cdot\mathbf{u}^h \right)$
 in (\ref{sec:Towards correct energy, GLSDSD}), where $\mathbf{e}_i$ is the $i$th Cartesian basis vector. 
 Using $\nabla \mathbf{e}_i=\B{0}$ and the pointwise divergence-free velocity (\ref{div free}), all diffusive 
 and pressure terms drop out and we are left with:
\begin{align}\label{global line ar mom conv term initial}
  \left(\mathbf{e}_i, \pd_t \mathbf{u}^h  + \pd_t \mathbf{u}' \right) 
  +\tfrac{1}{2} \left(\mathbf{e}_i, \left(\left(\mathbf{u}^h  
  + \mathbf{u}'\right)\cdot \nabla\right) \mathbf{u}^h\right) 
  + \left( \nabla \left(- \tfrac{1}{2}\mathbf{e}_i\cdot \mathbf{u}^h\right) ,  \mathbf{u}' \right)=  (\mathbf{e}_i,f).
\end{align}
Consider the convective term in isolation and write
\begin{align}\label{global line ar mom conv term}
  \tfrac{1}{2}\left(\mathbf{e}_i, \left(\left(\mathbf{u}^h  + \mathbf{u}'\right)\cdot \nabla\right) \mathbf{u}^h\right)  
  &=  \tfrac{1}{2}\left(\mathbf{e}_i , \nabla \cdot \left( \left(\mathbf{u}^h  + \mathbf{u}'\right) \otimes \mathbf{u}^h\right) \right)
   - \tfrac{1}{2}\left(\mathbf{e}_i ,  \left(\nabla \cdot   \left(\mathbf{u}^h  + \mathbf{u}'\right) \right)  \mathbf{u}^h\right)   \nonumber \\
 &= -  \tfrac{1}{2}\left(\nabla \mathbf{e}_i ,   \left(\mathbf{u}^h  + \mathbf{u}'\right) \otimes \mathbf{u}^h\right) 
  -\tfrac{1}{2} \left(\mathbf{e}_i ,  \left(\nabla \cdot   \left(\mathbf{u}^h  + \mathbf{u}'\right) \right)  \mathbf{u}^h\right)\nonumber \\
  &=    -  \tfrac{1}{2}\left(\mathbf{e}_i\cdot \mathbf{u}^h,  \nabla \cdot   \mathbf{u}' \right) \nonumber \\
 &=   \left( \nabla \left(\tfrac{1}{2}\mathbf{e}_i\cdot \mathbf{u}^h\right) ,  \mathbf{u}' \right).
\end{align}
Combining (\ref{global line ar mom conv term initial}) and (\ref{global line ar mom conv term}) leads to the balance
\begin{align}
  (\mathbf{e}_i,  \pd_t \mathbf{u}^h  + \pd_t \mathbf{u}') =  (\mathbf{e}_i,\mathbf{f}).
\end{align}
Linear momentum is thus conserved in terms of the total solution.
\subsection{Angular momentum}\label{subsec:Angular momentum}
Conservation of global angular momentum is a desirable property, certainly in rotating flows. 
It has been analyzed by Bazilevs et al. \cite{BaAk10} and Evans et al. \cite{Evans13unsteadyNS}. 
When using the appropriate weighting function spaces the formulation conserves angular momentum. 
The numerical results of Section \ref{sec:ns_case} are however not computed with these weighting 
function spaces. The demonstration of conservation of angular momentum follows the same ideas 
as \cite{BaAk10}. We set the weighting functions 
$\left(\mathbf{w}^h,q^h, \theta^h \right)=\left(\mathbf{x} \times \mathbf{e}_j,0,-\tfrac{1}{2}\left(\mathbf{x} \times \mathbf{e}_j\right)\cdot\mathbf{u}^h\right)$. 
By construction the gradient of the weighting function leads to a skew-symmetric tensor \cite{BaAk10}. 
As a result the gradient tensor is orthogonal to any symmetric tensor. Consequently the 
divergence, which is the trace of the gradient, is zero.

Employing these weighting functions in the weak form we arrive at
\begin{align}
  &( \mathbf{x} \times \mathbf{e}_j, \pd_t \mathbf{u}^h  + \pd_t \mathbf{u}')  
  +\tfrac{1}{2} (\mathbf{x} \times \mathbf{e}_j, ((\mathbf{u}^h  + \mathbf{u}')\cdot\nabla) \mathbf{u}^h)_{\Omega} 
  - \tfrac{1}{2} (\left((\mathbf{u}^h  + \mathbf{u}')\cdot \nabla\right) \left( \mathbf{x} \times \mathbf{e}_j\right), \mathbf{u}^h )\nonumber \\
  &- \left(\left((\mathbf{u}^h+\mathbf{u}')\cdot \nabla\right) \left(\mathbf{x} \times \mathbf{e}_j\right),\mathbf{u}'\right)
  -\tfrac{1}{2}\left(\nabla \left(\left( \mathbf{x} \times \mathbf{e}_j\right)\cdot \mathbf{u}^h\right),\mathbf{u}'\right)
  =  (\left(\mathbf{x} \times \mathbf{e}_j\right),\mathbf{f}).
\end{align}
Consider again the convective terms in isolation. \added[id=Rev.1]{Switching back to a conservative form, see (\ref{convective derivation}), yields an incompressibility term:}
\begin{align}
  &\tfrac{1}{2} (\mathbf{x} \times \mathbf{e}_j, ((\mathbf{u}^h  + \mathbf{u}')\cdot\nabla) \mathbf{u}^h) 
  - \tfrac{1}{2} (\left((\mathbf{u}^h  + \mathbf{u}')\cdot \nabla\right) \left( \mathbf{x} \times \mathbf{e}_j\right), \mathbf{u}^h )\nonumber \\
  &- \left(\left((\mathbf{u}^h+\mathbf{u}')\cdot \nabla\right) \left(\mathbf{x} \times \mathbf{e}_j\right),\mathbf{u}'\right)\nonumber\\
  =&-(\nabla \left(\mathbf{x} \times \mathbf{e}_j\right), (\mathbf{u}^h+\mathbf{u}')\otimes (\mathbf{u}^h+\mathbf{u}'))
  -\tfrac{1}{2}\left(\mathbf{u}^h,\left(\mathbf{x} \times \mathbf{e}_j\right)\nabla \cdot \left(\mathbf{u}^h+\mathbf{u}'\right)\right)\nonumber \\
    =&-(\nabla \left(\mathbf{x} \times \mathbf{e}_j\right), (\mathbf{u}^h+\mathbf{u}')\otimes (\mathbf{u}^h+\mathbf{u}')) 
    + \tfrac{1}{2}\left(\nabla \left( \left(\mathbf{x} \times \mathbf{e}_j\right)\cdot\mathbf{u}^h\right),\mathbf{u}'\right).
\end{align}
The antisymmetric tensor and the symmetric tensor in the first and second argument, \added[id=Authors]{respectively, 
cause} the first term to vanish. The \added[id=Rev.1]{incompressibility term} cancels with the choice of $\theta^h$ 
and the conservation of angular momentum is what remains:
\begin{align}
  ( \mathbf{x} \times \mathbf{e}_j, \pd_t \mathbf{u}^h  + \pd_t \mathbf{u}')  =(\mathbf{x} \times \mathbf{e}_j,\mathbf{f}).
\end{align}
\section{Numerical test case}
\label{sec:ns_case}
In this Section we test the GLSDD method (\ref{sec:Towards correct energy, GLSDSD}) on 
a three-dimensional Taylor--Green vortex flow at Reynolds number $Re=1600$. This test 
case is challenging and it is often employed to examine the performance of numerical 
algorithms for turbulence computations. It serves our purpose because (i) the energy 
behavior of a fully turbulent flow can be studied, (ii) reference data is available and 
(iii) the domain is periodic. Other boundary conditions than periodic ones are 
beyond the scope of this work.

\begin{figure}[h!]
    \begin{center}
    \begin{subfigure}[b]{0.48\textwidth}
         \includegraphics[scale=0.25]{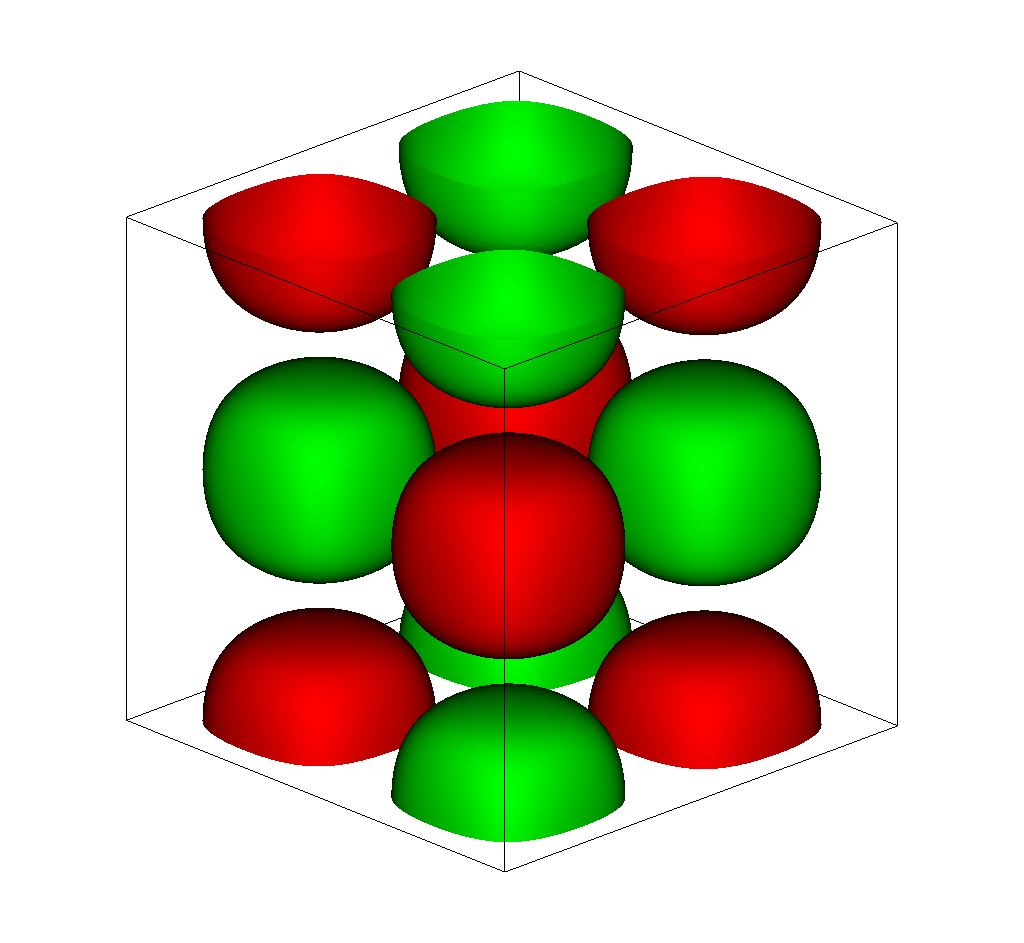}
        \caption{Laminar flow at $t=0$~s.}
    \end{subfigure}
    \begin{subfigure}[b]{0.48\textwidth}
        \includegraphics[scale=0.25]{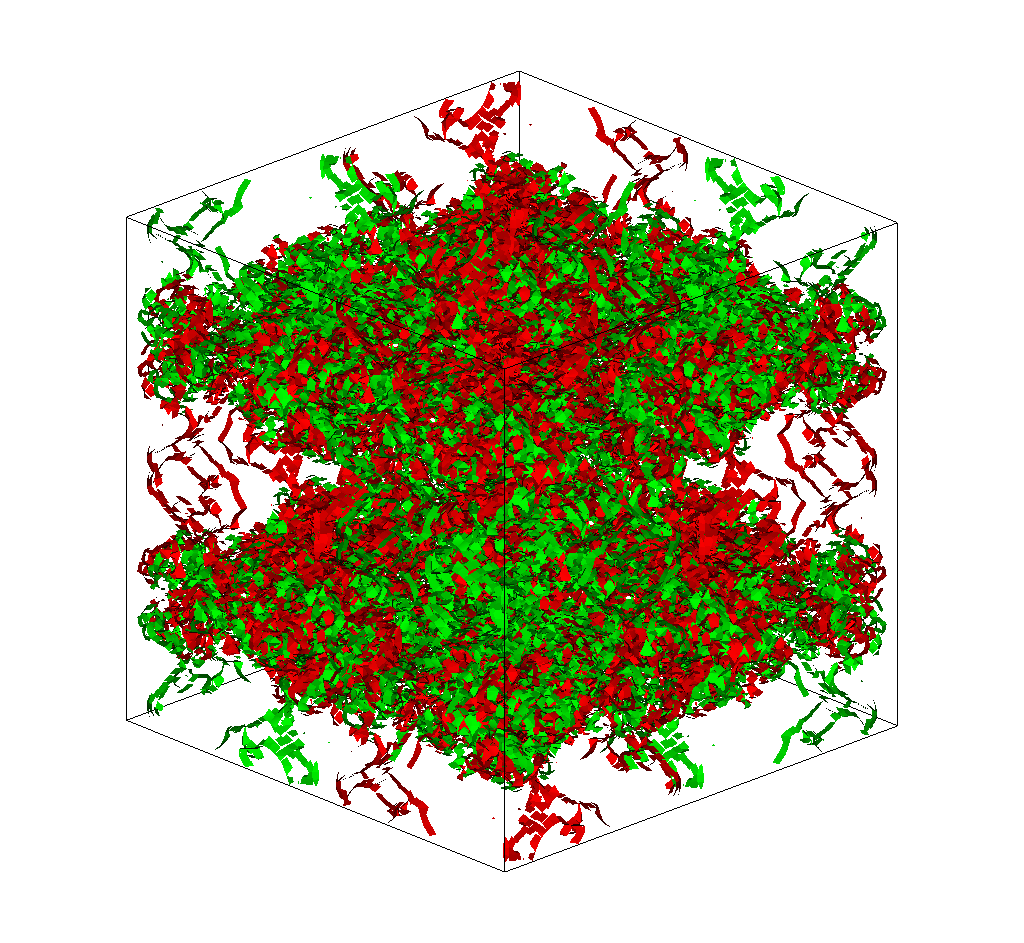}
        \caption{Fully turbulent flow at $t=10$~s.}
    \end{subfigure}
    \caption{Taylor--Green vortex flow at $Re=1600$. Iso-surfaces of z-vorticity.}\label{fig:TG vis}
    \end{center}
\end{figure}

The flow is initially of laminar type. As the time evolves, the vortices begin to 
evolve and roll-up. The vortical structures undergo changes and subsequently 
their structures breakdown and form distorted vorticity patches. The flow
transitions to one with a turbulence character; the vortex stretching causes 
the creation of small-scales.
The Taylor--Green vortex initial conditions are specified as follows:
\begin{subequations}
  \label{sec:Towards correct energy, GLS}
  \begin{alignat}{2}
    u(\mathbf{x},0)&=\sin(x)\cos(y)\cos(z),\\
    v(\mathbf{x},0)&= -\cos(x)\sin(y)\cos(z),\\
    w(\mathbf{x},0)&=0,\\
    p(\mathbf{x},0)&= \tfrac{1}{16}\left(\cos(2 x)+\cos(2 y)\right)\left(\cos(2z)+2\right).
  \end{alignat}
\end{subequations}
The physical domain is the cube $\Omega=\left[0,2\pi\right]^3$ with periodic 
boundary conditions. For this test case the viscosity is given by $\nu=\frac{1}{Re}$. 
Here we consider the transition phase for times $t \leq 10$~s. Figure \ref{fig:TG vis} 
shows the iso-surfaces of the z-vorticity of the initial condition (laminar flow) 
and the final configuration (fully turbulent flow).\\
\indent Due to the symmetric behavior of the flow, we are allowed to simulate only 
an eighth part of the domain. Hence, we take as computational domain 
$\Omega^h=\left[0,\pi\right]^3$ and apply no-penetration boundary conditions.
All the implementations employ NURBS basis functions that are mostly $C^1$-quadratic, 
however every velocity space is enriched to be cubic $C^2$ in the \added[id=Rev.2]{associated 
direction} \cite{Evans13steadyNS, Evans13unsteadyNS, buffa2011isogeometric, buffa2011isogeometric3D}. 
Note that conservation of angular momentum cannot be guaranteed, since the 
\added[id=Rev.1]{choice of the weighting function $\theta^h$ in section \ref{subsec:Angular momentum}} is not valid. We apply a standard
$L_2$-projection to set the initial condition on the mesh. For the time-integration
we employ the generalized-$\alpha$ method with the parameter choices of \cite{EAk17} 
which yield correct energy evolution. This method is stable and shows second-order temporal accuracy.
\added[id=Rev.1]{The resulting system of equations is solved with the standard flexible GMRES method with additive Schwartz preconditioning provided by Petsc \cite{PETSC,petsc-efficient}.}

We perform simulations with three different methods: (i) the classical Galerkin method,
(ii) the VMS method with static small-scales (VMSS), comparable with \cite{BaCaCoHu07} and 
(iii) the novel Galerkin/least-squares formulation with dynamic and divergence-free 
small-scales (GLSDD), i.e. form (\ref{sec:Towards correct energy, GLSDSD}). 
The DNS results of Brachet et al. \cite{brachet1983small} obtained with a spectral 
method on a fine $512^3$-mesh serve as reference data (ref). 

\begin{figure}[h!]
    \begin{center}
    \begin{subfigure}[b]{0.49\textwidth}
  \includegraphics[scale=0.625]{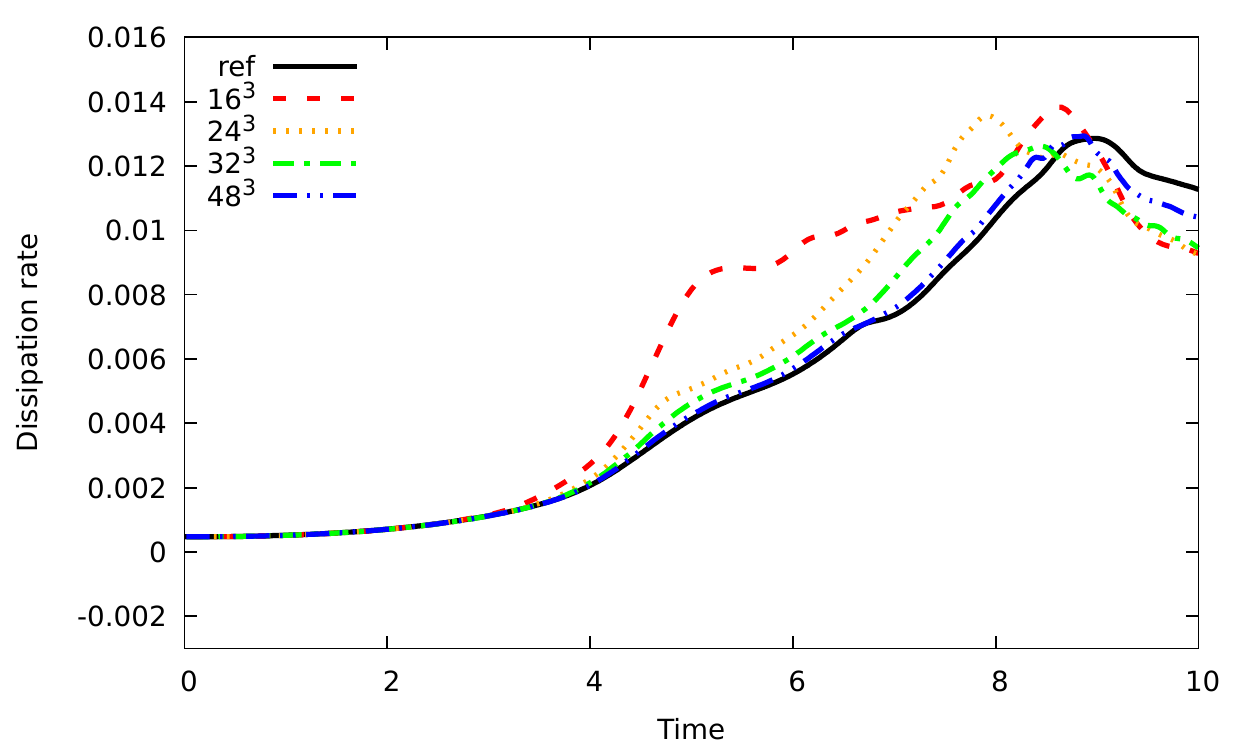}
  \caption{Dissipation rate}
  \label{fig: TG dissipation mesh ref} 
    \end{subfigure}
    \begin{subfigure}[b]{0.49\textwidth}
      \includegraphics[scale=0.625]{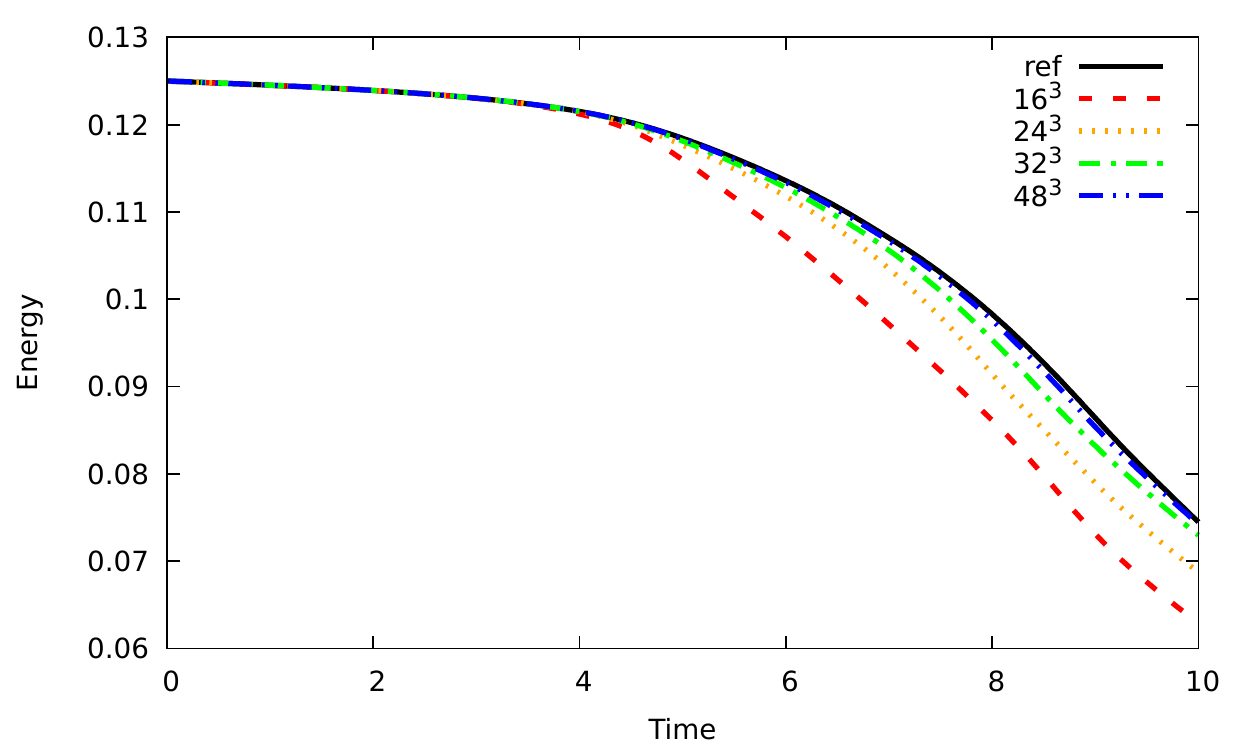}
      \caption{Energy evolution}
  \label{fig: TG energy mesh ref}
    \end{subfigure}
    \caption{Taylor--Green vortex flow at $Re=1600$ mesh convergence. The 
                 GLS method with dynamic divergence-free small-scales.}
    \label{fig: TG mesh ref}
    \end{center}
\end{figure}

First, we perform a brief mesh refinement study for the novel method. 
Figure \ref{fig: TG mesh ref} shows  mesh refined results for the novel 
GLSDD method (\ref{sec:Towards correct energy, GLSDSD}). 
For this purpose meshes with $16^3$, $24^3$, $32^3$, $48^3$ \added[id=Authors]{elements} have been employed. Clearly, the energy behavior on the coarsest
two meshes is quite off. The finer meshes are able to closely capture the
turbulence character of the flow. In the following we therefore use 
meshes of $32^3$ or $48^3$ elements.

\begin{figure}[h!]
    \begin{center}
    \begin{subfigure}[b]{0.49\textwidth}
  \includegraphics[scale=0.625]{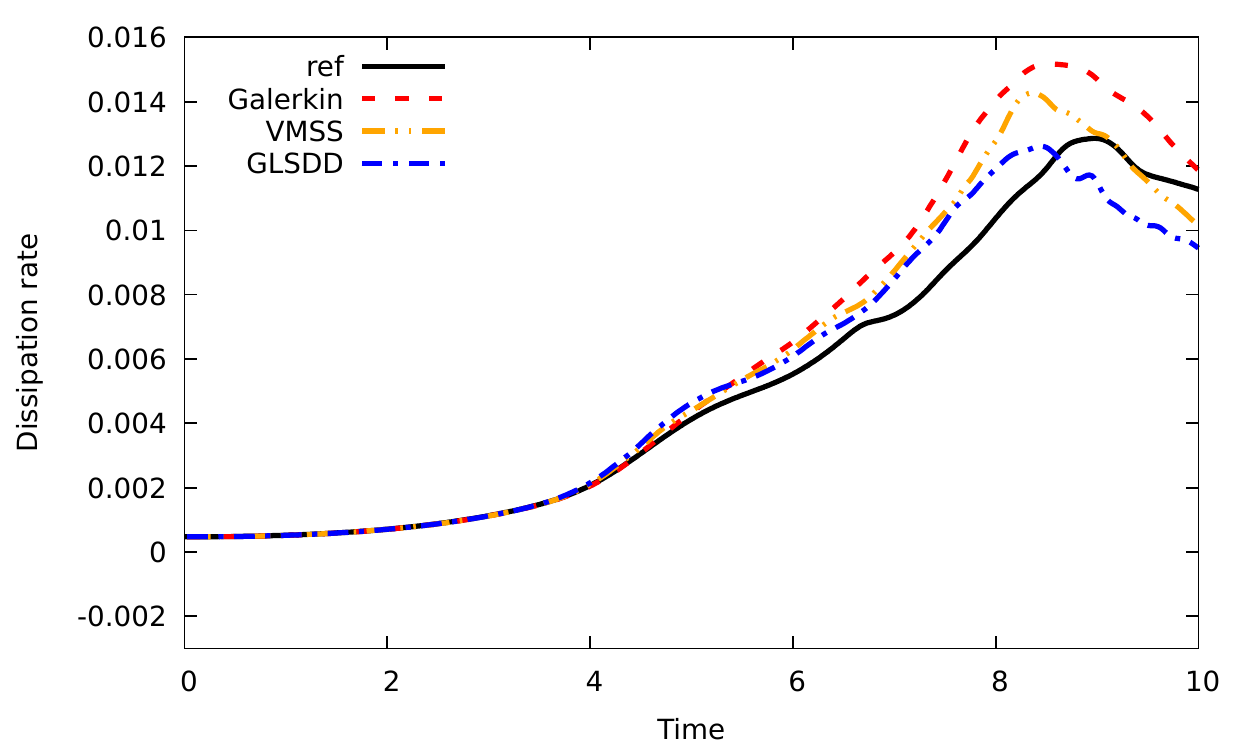}
  \caption{Dissipation rate}
  \label{fig: TG dissipation 32}
    \end{subfigure}
    \begin{subfigure}[b]{0.49\textwidth}
      \includegraphics[scale=0.625]{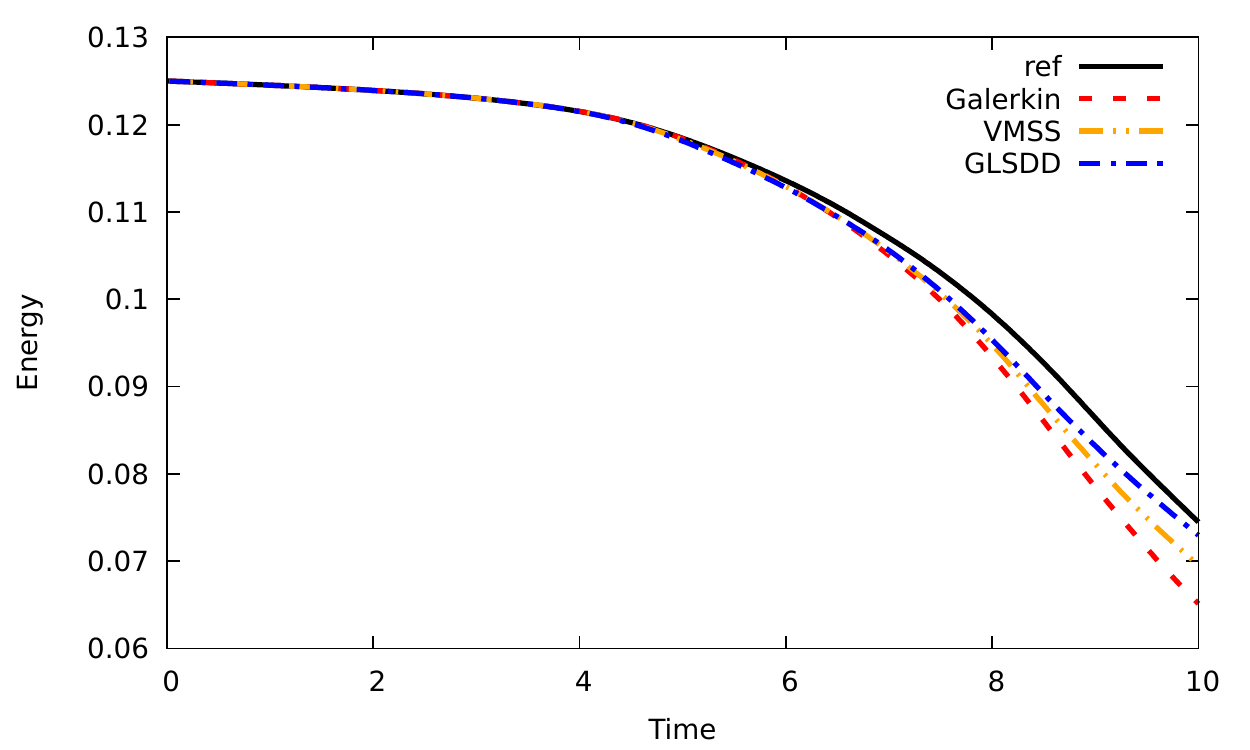}
      \caption{Energy evolution}
      \label{fig: TG energy 32}
    \end{subfigure}
    \caption{Taylor--Green vortex flow at $Re=1600$ on $32^3$-mesh for 
    various methods: the Galerkin method, the VMS method with static 
    small-scales and the GLS method with dynamic divergence-free small-scales.}\label{fig: TG 32}
    \end{center}
\end{figure}

\begin{figure}[h!]
    \begin{center}
    \begin{subfigure}[b]{0.49\textwidth}
  \includegraphics[scale=0.625]{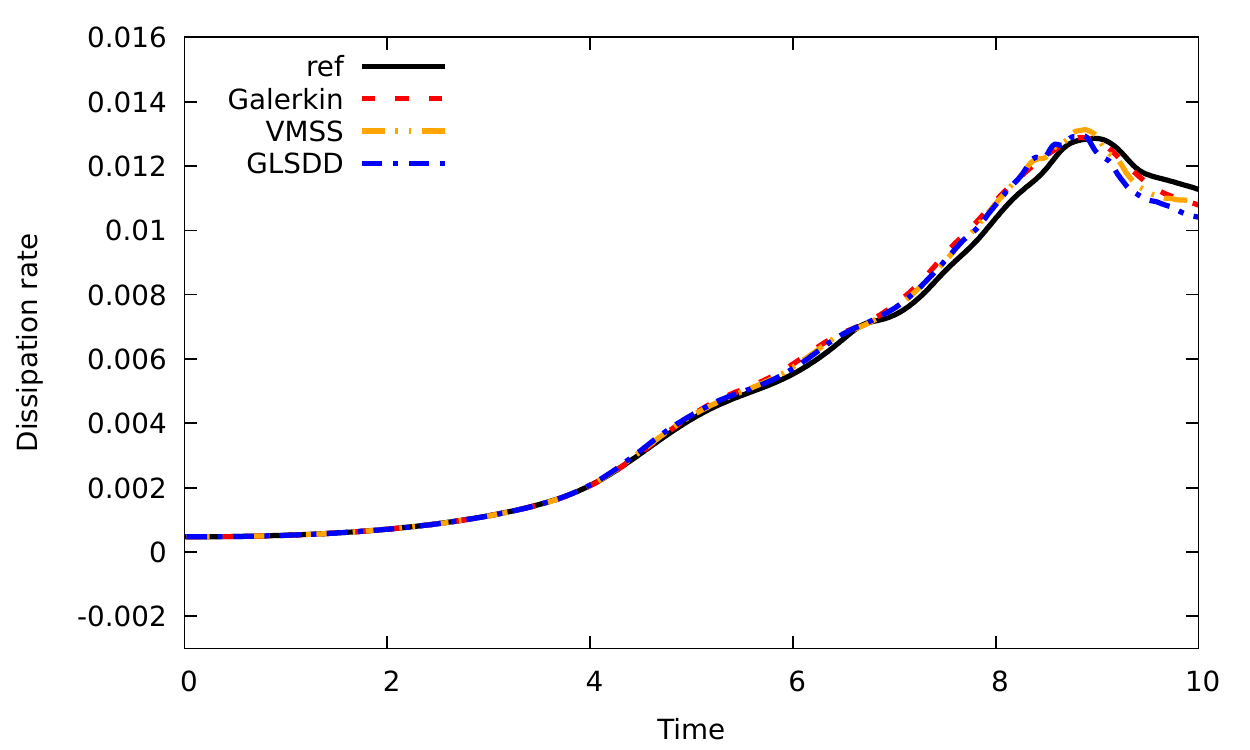}
  \caption{Dissipation rate}
  \label{fig: TG dissipation 48}
    \end{subfigure}
    \begin{subfigure}[b]{0.49\textwidth}
      \includegraphics[scale=0.625]{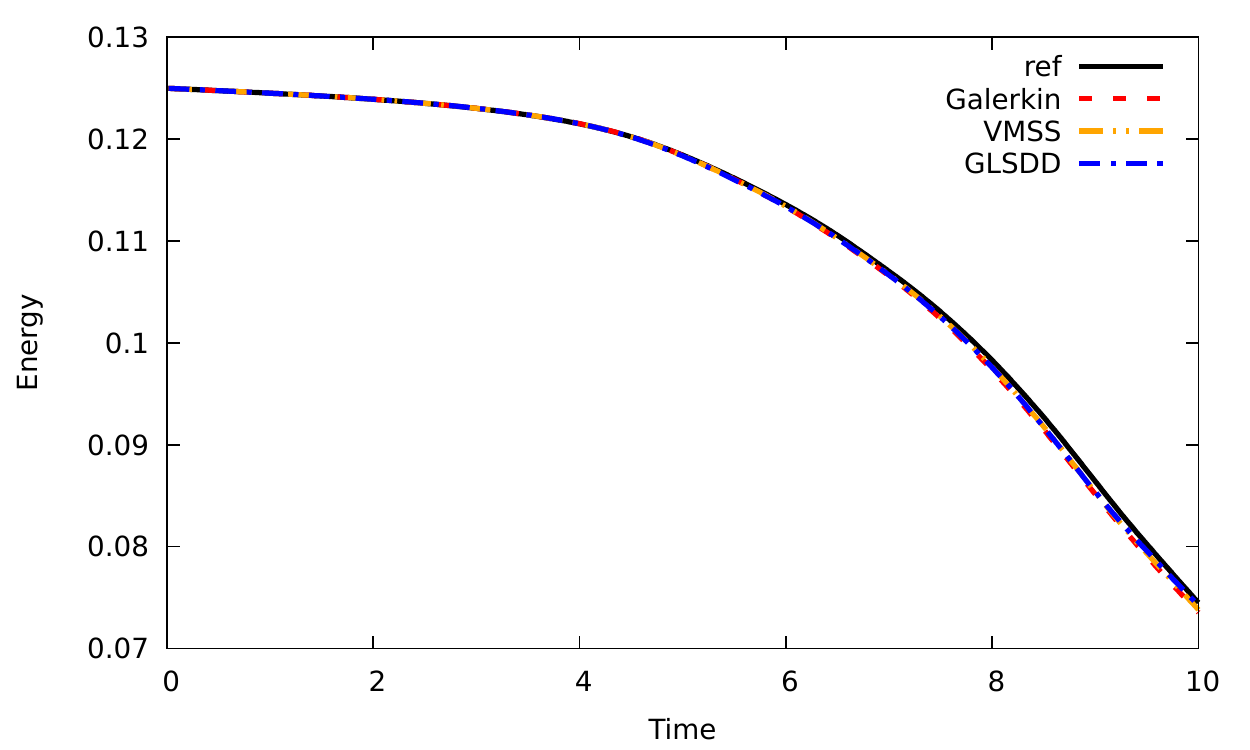}
      \caption{Energy evolution}
      \label{fig: TG energy 48}
    \end{subfigure}
    \caption{Taylor--Green vortex flow at $Re=1600$ on $48^3$-mesh for various 
    methods: the Galerkin method, the VMS method with static small-scales and 
    the GLS method with dynamic divergence-free small-scales.}\label{fig: TG 48}
    \end{center}
\end{figure}

We compare the results of the novel GLSDD method with the VMSS and the
Galerkin approach. The simulations are carried out on a mesh of $32^3$ 
elements, i.e. the mesh size is $h=\frac{\pi}{32}$, and on a slightly finer mesh 
of $48^3$ elements. The time-step is taken as $\Delta t = \frac{ 4 h}{5\pi}$, 
i.e. the initial CFL-number is roughly $0.25$. In the Figures \ref{fig: TG 32}-\ref{fig: TG 48} 
we visualize the time history of the kinetic energy and kinetic energy dissipation 
rate for each of the three methods and the reference data.

The Figure \ref{fig: TG 32} shows that each of the methods is able to roughly 
capture the energy behavior on the coarse mesh. The dissipation peek 
appears too early in time for each of the simulations. The Galerkin method 
displays the least accurate results, it overpredicts the dissipation rate. 
The VMSS method performs a bit better at all times. The novel GLSDD 
approach demonstrates an even closer agreement with the reference results. 
The results on the finer mesh, in Figure \ref{fig: TG 48}, reveal 
almost no difference with the reference data.

\begin{figure}[h!]
    \begin{center}
    \begin{subfigure}[b]{0.49\textwidth}
  \includegraphics[scale=0.625]{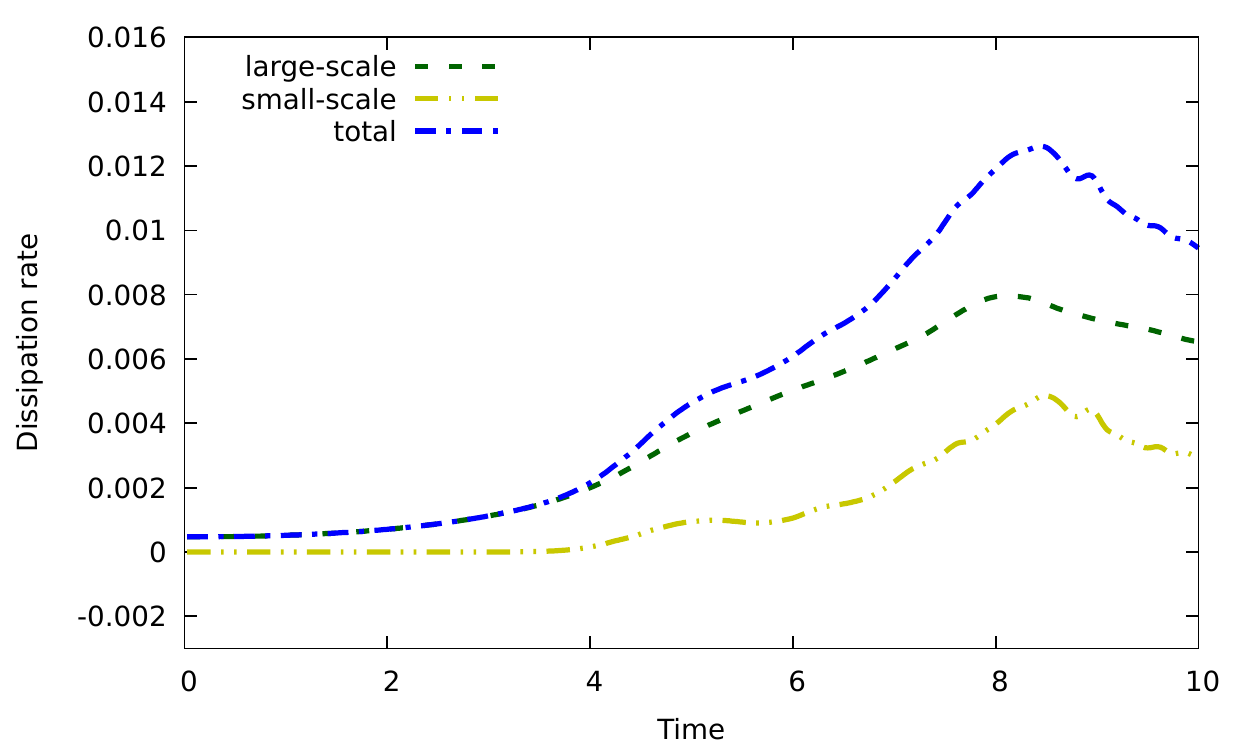}
  \caption{Large- and small-scale contributions}
  \label{fig: TG dissipation 32 split}
    \end{subfigure}
    \begin{subfigure}[b]{0.49\textwidth}
      \includegraphics[scale=0.625]{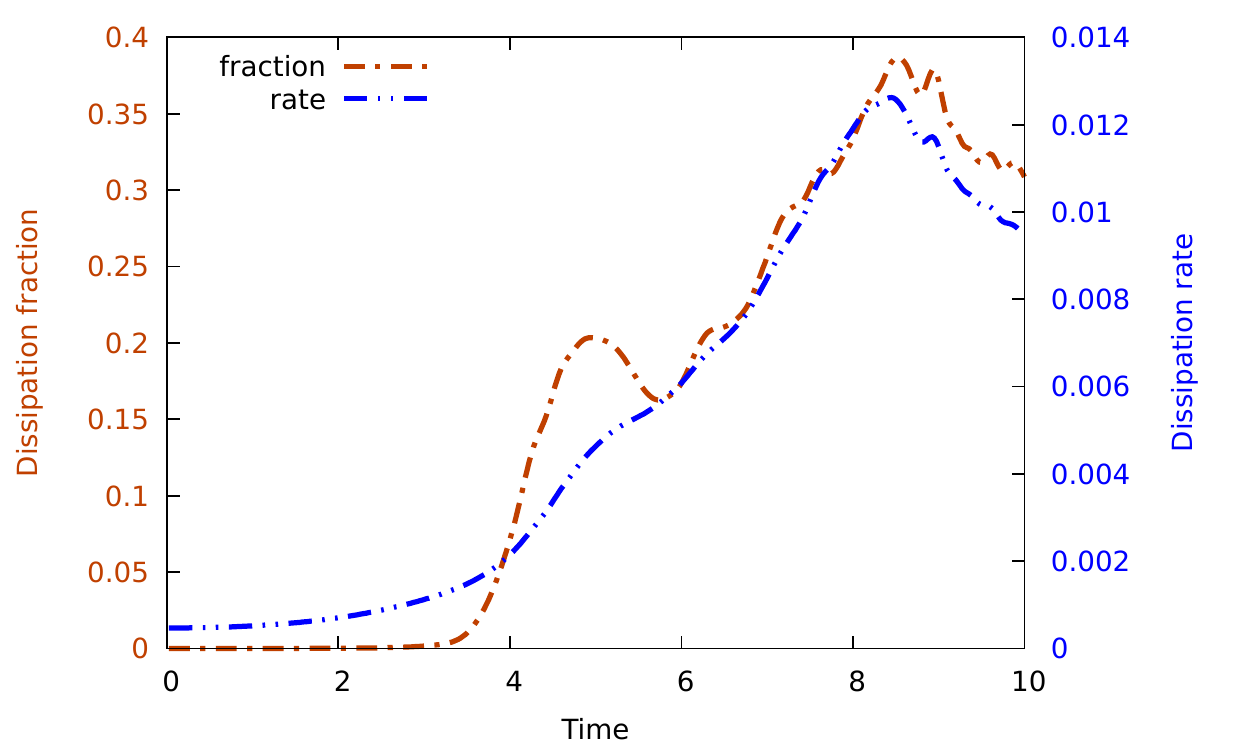}
      \caption{The small-scale dissipation fraction}
      \label{fig: TG energy 32 rate fraction}
    \end{subfigure}
    \caption{Taylor--Green vortex flow at $Re=1600$ on $32^3$-mesh with the 
    GLSDD method: (a) large-scale and small-scale parts of the dissipation 
    rate and (b) their fraction.}\label{fig: TG 32 split fraction}
    \end{center}
\end{figure}

In the following we further analyze the contributions of the dissipation
 rate (\added[id=Rev.1]{on the course mesh}). The dissipation rate of the Galerkin 
 method only consists of the large-scale/physical dissipation 
 $\|\nu^{1/2} \nabla \mathbf{u}^h\|_{\Omega}^2$. In contrast, 
 the dissipation of the GLSDD method is composed of a large-scale and a small-scale contribution:
\begin{align}
 \dfrac{{\rm d}}{{\rm d}t} E_\Omega^{\text{GLSDD}}= & - \|\nu^{1/2}\nabla \mathbf{u}^h\|_{\Omega}^2-\|\tau_M^{-1/2}\mathbf{u}'\|_{\tilde{\Omega}}^2.
\end{align}
In Figure \ref{fig: TG 32 split fraction} we display the temporal evolution of 
both parts and the small-scale dissipation fraction \newline 
$(\|\tau_M^{-1/2}\mathbf{u}'\|_{\tilde{\Omega}}^2)/(\|\nu^{1/2}\nabla \mathbf{u}^h\|_{\Omega}^2+\|\tau_M^{-1/2}\mathbf{u}'\|_{\tilde{\Omega}}^2)$. 
In the laminar regime ($t<3$) the small-scale contribution is negligible. 
When the flow has a more turbulent character the contribution of the 
small-scales is substantial: the maximum of the dissipation fraction exceeds $0.35$.

Lastly, we focus on the energy dissipation of the VMSS formulation. 
The derivation follows the same steps used throughout this paper. 
\added[id=Rev.1]{One might argue that the energy could also be solely based on the large-scales. This is what we do here.} Its evolution reads:
\begin{align}\label{energy evolution stat VMS}
 \dfrac{{\rm d}}{{\rm d}t} E_\Omega^{h,\text{VMSS}}=& 
 - \|\nu^{1/2}\nabla \mathbf{u}^h\|_{\Omega}^2 
 -\|\tau_M^{-1/2}\mathbf{u}'\|_{\tilde{\Omega}}^2  
 + (\nu \Delta \mathbf{u}^h,\mathbf{u}')_{\tilde{\Omega}} 
 -\left( \mathbf{u}',\partial_t \mathbf{u}^h \right)_{\tilde{\Omega}}\nonumber \\
                     &+(\nabla \cdot \mathbf{u}^h,p')_{\tilde{\Omega}} 
 +\left(\nabla \mathbf{u}^h , (\mathbf{u}^h + \mathbf{u}') \otimes (\mathbf{u}^h + \mathbf{u}') \right)_{\tilde{\Omega}} 
 - \left(\mathbf{u}',(\mathbf{u}^h+\mathbf{u}')\cdot \nabla \mathbf{u}^h \right)_{\tilde{\Omega}}.
\end{align}

\begin{figure}[h!]
    \begin{center}
    \begin{subfigure}[b]{0.49\textwidth}
  \includegraphics[scale=0.625]{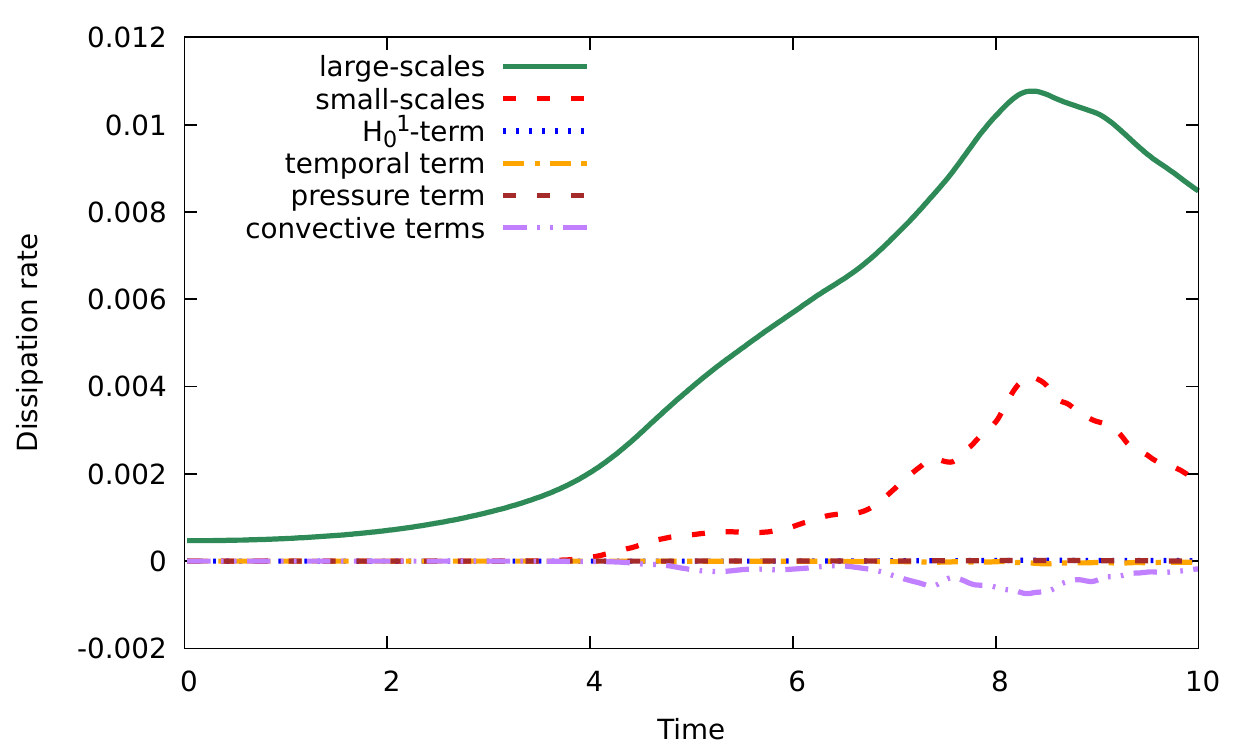}
  \caption{All terms}
  \label{fig: TG dissipation 32 split}
    \end{subfigure}
    \begin{subfigure}[b]{0.49\textwidth}
      \includegraphics[scale=0.625]{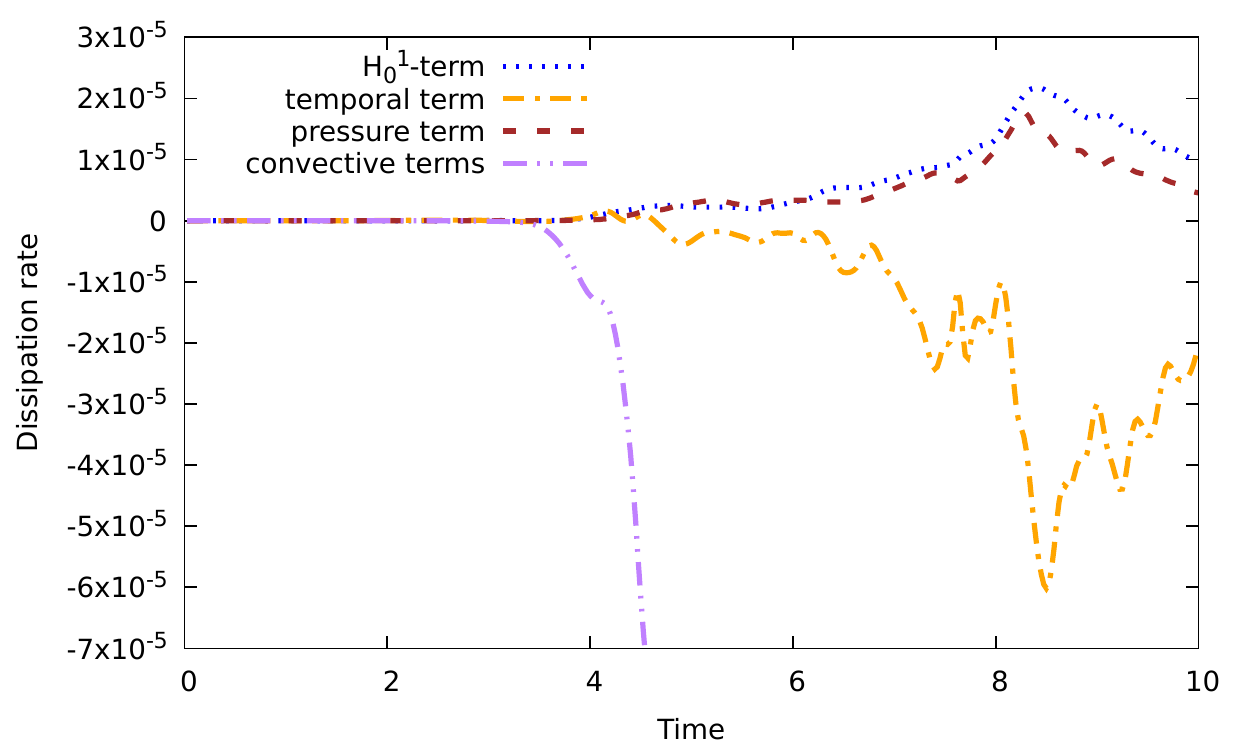}
      \caption{Zoom unwanted terms}
      \label{fig: TG energy 32 rate fraction}
    \end{subfigure}
    \caption{Taylor--Green vortex flow at $Re=1600$ on $32^3$-mesh with the 
    VMSS method: energy dissipation of separate terms.}\label{fig: TG VMSS32 smallscaledissipation}
    \end{center}
\end{figure}

Figure \ref{fig: TG VMSS32 smallscaledissipation} shows the contribution of the separate terms. 
The two desired dissipation terms are clearly dominant. The small-scale 
dissipation is smaller than the large-scale dissipation, however it has a 
significant contribution. Although the contributions are small, the unwanted 
terms can create artificial energy.

\section{Conclusions}
\label{sec:CONC}
We continued the study initiated in \cite{EAk17} \added[id=Rev.1]{concerning} the construction of 
methods displaying correct-energy behavior. In this paper we have 
applied the developed methodology to the incompressible Navier--Stokes 
equations. It clearly shows \added[id=Rev.1]{that the link between the methods VMS, SUPG and GLS, established in \cite{EAk17}, is also valid for the incompressible Navier--Stokes equations.}

The novel GLSDD methodology employs divergence-conforming NURBS 
basis functions and uses a Lagrange multiplier setting to ensure divergence-free 
small-scales. Furthermore, it enjoys the favorable behavior of the dynamic 
small-scales and reduces to the Galerkin method in the Stokes regime. 
These properties all emerge from the correct-energy design condition. 
A pleasant byproduct of the method is the conservation of linear momentum. 
The conservation of angular momentum can be achieved when employing 
the appropriate weighting function spaces. The numerical results support the 
theoretical framework in that the energy behavior improves upon the VMS 
method with static small-scales. The variational multiscale method with
 static small-scales has unwanted small-scale contributions 
 which create artificial energy.

The novel formulation requires a bit more effort to \added[id=Rev.1]{implement} compared 
\added[id=Rev.1]{to} the variational multiscale method with static small-scales. One has to 
include an additional variable to ensure the divergence-free behavior of the 
small-scales. In addition the formulation \added[id=Rev.1]{needs} to be equipped with the 
dynamic small-scale model. However, the resulting system of equations
 does not demand a sophisticated preconditioner; we have employed the
  standard ASM (Additive Schwarz Method) technique. In our opinion, the 
  accuracy gain outweighs the little \added[id=Rev.1]{extra implementation 
  effort} and calculation cost.

We see several directions for future work. The first concerns the 
development of a method displaying correct energy behavior at the 
boundary, in particular \added[id=Rev.1]{when} using the weak imposition of Dirichlet 
boundary conditions. This allows to test the effect of correct energy 
behavior on wall-bounded turbulent flow problems. Another extension 
is correct energy behavior for free-surface flow computations. This is an 
important step, since artificial energy creation can yield highly instable 
behavior, as demonstrated in \cite{AkBaBeFaKe12}. We have work on both extensions in 
progress and aim to report on it in the near future.

\section*{Acknowledgment}
\label{sec:ack}

The authors are grateful to the Delft University of Technology
for its support.

\appendix

\section{Galerkin/least-squares formulation with dynamic divergence-free small-scales}
\label{Appendix: Galerkin/least-squares formulation with dynamic divergence-free small-scales}
We repeat the Galerkin/least-squares formulation with dynamic divergence-free 
small-scales (GLSDD), i.e. form  (\ref{sec:Towards correct energy, GLSDSD}), to 
provide an overview of the separate terms. The formulation is of skew-symmetric type, 
applies GLS stabilization and uses divergence-free dynamic small-scales. The method requires a stable velocity--pressure pair and reads:\\

\textit{Find $\left(\mathbf{u}^h,p^h, \zeta^h\right) \in \WW^h\times \mathcal{P}^h$ 
such that for all $\left(\mathbf{w}^h,q^h,\theta^h \right) \in \WW^h\times\mathcal{P}^h$,}
\begin{subequations}
\label{eq:omegai}
\begin{align}
\label{sec:Towards correct energy, GLSDSD2}
\begin{split}
\left( \mathbf{w}^h, \pd_t  \mathbf{u}^h  \right)_{\Omega} +{}& \left(\mathbf{w}^h, \pd_t \mathbf{u}'\right)_{\tilde{\Omega}}\\[6pt]
 +& \tfrac{1}{2}(\mathbf{w}^h, (\mathbf{u}^h+\mathbf{u}') \cdot \nabla \mathbf{u}^h)_{\Omega}
 -\tfrac{1}{2}((\mathbf{u}^h+\mathbf{u}')\cdot\nabla\mathbf{w}^h, \mathbf{u}^h)_{\Omega}
- \left( \left(\mathbf{u}^h  + \mathbf{u}'\right)\cdot \nabla \mathbf{w}^h  ,  \mathbf{u}' \right)_{\tilde{\Omega}}\\[6pt]
+&(\nabla \mathbf{w}^h, 2 \nu \nabla^s \mathbf{u}^h )_{\Omega} 
+ \left( \nu \Delta \mathbf{w}^h, \mathbf{u}'\right)_{\tilde{\Omega}} \\[6pt]
+&(q^h,\nabla \cdot \mathbf{u}^h)_{\Omega}- (\nabla \cdot \mathbf{w}^h, p^h )_{\Omega}  
+\left( \nabla \theta^h, \mathbf{u}'\right)_{\tilde{\Omega}}=(\mathbf{w},\mathbf{f})_{\Omega},
\end{split}
\\[6pt]
\label{eq:omega2}
\begin{split}
\pd_t \mathbf{u}'+\tau_M^{-1} \mathbf{u}' + \nabla \zeta^h  {}&+ \B{r}_M =0,
\end{split}
\end{align}
\end{subequations}
where momentum residual is
\begin{equation}
 \B{r}_M= \pd_t \mathbf{u}^h + \left(\left(\mathbf{u}^h+\mathbf{u}'\right)\cdot\nabla\right) \mathbf{u}^h 
 + \nabla p^h - \nu \Delta \mathbf{u}^h - \mathbf{f}.
\end{equation}
The separate terms of (\ref{sec:Towards correct energy, GLSDSD2}) are from 
left to right: the temporal terms, the convective contributions, the viscous contributions, 
the incompressibility constraint, the pressure term, the divergence-free small-scale 
velocity constraint and the forcing term. This form follows the correct-energy evolution (on a local scale):
\begin{align}\label{energy evolution GLSDD Appendix}
 \dfrac{{\rm d}}{{\rm d}t} E_\omega=& - \|\nu^{1/2}\nabla \mathbf{u}^h\|_{\omega}^2+
  (\mathbf{u}^h,\mathbf{f})_{\omega}- (1,F_\omega^h)_{\chi_\omega}\nonumber  \\
                     & -\|\tau_M^{-1/2}\mathbf{u}'\|_{\tilde{\omega}}^2 + (\mathbf{u}',\mathbf{f})_{\tilde{\omega}},
\end{align}
and possesses the conservation properties of Section \ref{sec:CP}.

\section{Definition dynamic stabilization parameter}
\label{Appendix: Definition stabilization parameters}
The dynamic stabilization parameter $\tau_{M}$ is the discrete 
approximation of the inverse of the convective and viscous parts 
of momentum Navier--Stokes operator. It mirrors the dynamic 
stabilization parameter of convection--diffusion equation (see \cite{EAk17}). 
The continuity stabilization parameter $\tau_C$ is on its turn the discrete 
approximation of the inverse of the divergence operator, here we use the 
objective definition introduced in \cite{BaAk10}. The parameters take the form:
\begin{subequations}
  \begin{alignat}{1} \label{eq:tau_static 2}
  \tau_{M}=&\left(\tau^{-2}_{\text{conv}}+\tau^{-2}_{\text{visc}}\right)^{-1/2},\\
  \tau_{C}=&\left(\tau_M \sqrt{ \mathbf{G}:\mathbf{G}}\right)^{-1},
\end{alignat}
\end{subequations}
where the convective and viscous contributions of $\tau_{M}$ are
\begin{subequations}
  \begin{alignat}{1}
  \tau_{\text{conv}}^{-2}=& 4 \mathbf{u}\cdot \mathbf{G} \mathbf{u},\\
  \tau_{\text{visc}}^{-2}=& C_I\nu^2 \mathbf{G}:\mathbf{G}.
\end{alignat}
\end{subequations}
Here the following definition is employed:
\begin{subequations}
  \begin{alignat}{1}
  \mathbf{G}&=\frac{\pd \boldsymbol{\xi}}{\pd \bx}^T\frac{\pd \boldsymbol{\xi}}{\pd \bx},\\
  \mathbf{G}:\mathbf{G} &= \displaystyle\sum_{i,j=1}^3 G_{ij} G_{ij},
\end{alignat}
\end{subequations}
where $\pd \boldsymbol{\xi}/\pd \bx$ is the inverse Jacobian of the map 
between the elements in the reference and physical domain. The positive 
constant $C_I$ is determined by an inverse estimate.

\bibliographystyle{unsrt}
\bibliography{references}

\end{document}